\documentclass[review, 12pt]{elsarticle}
\usepackage{lineno,hyperref,amsmath,amsthm}
\usepackage[noend]{algpseudocode}
\usepackage{algorithmicx,algorithm}
\usepackage{booktabs}
\usepackage{color}
\usepackage{graphics}
\usepackage{subfig}
\usepackage{amsfonts,amssymb}
\usepackage{bm}

\modulolinenumbers[1]

\journal{} 

\begin{document}
	
\begin{frontmatter}
	
\title{Structural stochastic responses determination via a sample-based stochastic finite element method}

\author[a,b]{Zhibao Zheng\corref{CorrespondingAuthor}}
\cortext[CorrespondingAuthor]{Corresponding author}
\ead{zhibaozheng@hit.edu.cn}


\address[a]{School of Civil Engineering, Harbin Institute of Technology, Harbin 150090, China}

\address[b]{Key Lab of Structures Dynamic Behavior and Control, Harbin Institute of Technology, Ministry of Education, Harbin 150090, China}

\begin{abstract} 
	This paper presents a new stochastic finite element method for computing structural stochastic responses. The method provides a new expansion of stochastic response and decouples the stochastic response into a combination of a series of deterministic responses with random variable coefficients. A dedicated iterative algorithm is proposed to determine the deterministic responses and corresponding random variable coefficients one by one. The algorithm computes the deterministic responses and corresponding random variable coefficients in their individual space and is insensitive to stochastic dimensions, thus it can be applied to high dimensional stochastic problems readily without extra difficulties. More importantly, the deterministic responses can be computed efficiently by use of existing Finite Element Method (FEM) solvers, thus the proposed method can be easy to embed into existing FEM structural analysis softwares. Three practical examples, including low-dimensional and high-dimensional stochastic problems, are given to demonstrate the accuracy and effectiveness of the proposed method.
\end{abstract}

				
\end{frontmatter}


\section{Introduction} \label{sec1}
Ordinary or partial differential equations (PDEs) are powerful tools to describe many real life engineering and scientific processes. A wide body of numerical methods based on finite differences, finite elements, and boundary elements are available to approximately solve the governing equations for the response quantities of interest. In particular, Finite element method (FEM) has become state-of-the-art since it offers a simple way to solve very high resolution models in various computational physics problems, ranging from structural mechanics, thermodynamics to nano-bio mechanics~\citep{hughes2012finite}. Nevertheless, the FEM is deterministic in nature and is therefore limited to describe the general characteristics of a real life system. The considerable influence of inherent uncertainties on system behavior has led the scientific community to recognize the importance of a stochastic approach to realistic engineering systems~\citep{Dai2015A}. More than ever, the goal becomes to represent and propagate uncertainties from the available data to the desired results through PDEs within the framework of stochastic equations~\citep{najm2009uncertainty, le2010spectral}.

The modelling of uncertainties consists in defining a suitable probability space $( {\Theta, \Sigma, \cal{P}} )$, where $\Theta$ denotes the space of elementary events, $\Sigma$ is a $\sigma $-algebra defined on $\Theta$ and $\cal{P}$ is a probability measure. In this paper, we consider the structural stochastic response $u\left( \theta  \right)$ of the problem is a stochastic function, with value in a certain function space, which has to verify almost surely the stochastic partial differential equations (SPDEs) discretized by stochastic finite element equations~\citep{Xiu2002The} as
\begin{equation} \label{eq:SFEM}
  K\left( \theta \right)u\left( \theta  \right) = F\left( \theta \right)
\end{equation}
where $K\left( \theta \right)$ is an operator representing properties of the physical model under investigation, which can be considered the stochastic stiffness matrix and $F\left( \theta \right)$ is a stochastic load vector. Randomness on the model can be formalized as a dependency of the operator and loads on the elementary event $\theta  \in \Theta$ and it's a great challenge to solve Equation \eqref{eq:SFEM} in the high dimensional stochastic space $\Theta$.

As an extension of deterministic FEM, stochastic finite element method (SFEM)~\citep{ghanem2003stochastic, Matthies2005Galerkin} has become a common tool for the solution of Equation \eqref{eq:SFEM}. Given the representation of uncertain system parameters and environmental source in terms of random fields, it becomes possible to integrate discretization methods for the response and random fields to arrive at a system of random algebraic equations. Two prominent variants of the SFEM are the non-intrusive methods and the Galerkin-type methods. Although various non-intrustive methods, e.g., Monte Carlo simulation~\citep{Papadrakakis1996Robust}, or regression and projection methods~\citep{blatman2008sparse}, can be readily applied to compute the response statistics to an arbitrary degree of accuracy, this is the method of last resort since the attendant computational cost can be prohibitive for real life problems. 

The Galerkin-type spectral methods~\citep{Xiu2002The, ghanem2003stochastic, Nouy2009Recent}, which are developed for linear SFEM, provide an explicit functional relationship between the random input and output, hence allow easy evaluation of the statistics of the stochastic system response. These methods transform Equation \eqref{eq:SFEM} arising from spatial discretization of SPDEs into a deterministic finite element equation by stochastic Galerkin projection, but the size of the deterministic finite element equation is significantly higher than that of the original SPDEs. Although several iterative solvers have been developed to decrease the substantial computational requisite~\citep{Ghanem1996Precondition, keese2005hierarchical}, the difficulty to build efficient preconditioners and memory requirements still limit their use to small-scale and low-dimensional stochastic problems. Also, the Curse of Dimensionality in stochastic spaces makes these methods more inefficient. For this line of approach to be successful in practice, it is crucial to have general-purpose and highly efficient numerical schemes for the solution of stochastic finite element equation \eqref{eq:SFEM}.

In this article, we develop a highly efficient numerical method for the explicit and high precision solution of Equation \eqref{eq:SFEM} with application to structural responses that involve uncertainties. An universal construct of solution to stochastic finite element (SFE) equations is firstly developed, which is independent on the types of SFE problems. Based on the construct of this solution, we further develop a numerical algorithm for solving SFE equations. The representations of the stochastic solutions are applicable for high-dimensional stochastic problems, and more importantly, the stochastic analysis and deterministic analysis in the solution procedure can thus be implemented in individual space. In this way, the proposed algorithm for the solution of SFE equation integrate the advantages of the non-instrusive methods and the Galerkin-type methods simultaneously, and thus have great potential for uncertainty quantifications in structural analysis.

The paper is organized as follows: Section \ref{sec2} briefly introduces the series expansion methods of random fields simulation and the derivation of stochastic finite element equations. A new method outlines for solving stochastic finite element equations is described in Section \ref{sec3}. Following this, the algorithm implementation of the proposed method is elaborated in Section \ref{sec4}. Three practical problems are used to demonstrate the accuracy and effectiveness of the proposed method in Section \ref{sec5}. Some conclusions and prospects are discussed in Sections \ref{sec6}.

\section{Stochastic finite element method} \label{sec2}
\subsection{Random fields expansion} \label{sec21}
In the framework of SFEM for structural analysis, uncertain physical parameters usually consists of the Young modulus, Poisson's ratio, yield stress, cross section geometry of physical systems, earthquake loading, wind loads, etc. In most cases, due to the lack of relevant experimental data, assumptions are made regarding probabilistic characteristics of random fields, sunch as Gaussian or non-Gaussian, stationary or non-stationary, etc~\citep{Dai2017A, Stefanou2009The}. The first step in applying the FEM to problems involving one or more of the random parameters is to model random fields based on the assumptions of probabilistic characteristics, thus random fields discretization is a key step in the numerical solutions of stochastic finite element equations. In order to derive stochastic finite element equations effectively, explicit expressions of random fields are also crucial. In general, we represent random fields by an enumerable set of random variables, and the series expansion of a second-order random field $\omega (x,\theta )$, which is indexed on a bounded domain $D$, can be expressed as
\begin{equation} \label{RF}
  \omega \left( {x,\theta } \right) = \sum\limits_{i = 0}^M {{\xi _i}\left( \theta  \right){\omega _i}\left( x \right)}
\end{equation}
where $\left\{ {{\xi _i}\left( \theta  \right)} \right\}_{i = 0}^M$ and $\left\{ {{\omega _i}\left( x \right)} \right\}_{i = 0}^M$ are random variables and deterministic functions, respectively, and $M$ is the number of retained items. Equation \eqref{RF} can be obtained by some methods for discretization of random fields. Various discretization techniques are available in the literature for approximating random fields including shape function methods, optimal linear estimation, weighted integral methods, orthogonal series expansion~\citep{sudret2000stochastic, 2010Identification}.

As a special case of the orthogonal series expansion, Karhunen-Lo\'{e}ve expansion is the most commonly used method in SFEM, and it has a form as 
\begin{equation} \label{KL}
  \omega \left( {x,\theta } \right) = {\omega _0}\left( x \right){\rm{ + }}\sum\limits_{i = 1}^M {{\xi _i}\left( \theta  \right)\sqrt {{\lambda _i}} {\omega _i}\left( x \right)}
\end{equation}
where ${\omega _0}\left( x \right)$ is the mean function of the random field $\omega \left( {x,\theta } \right)$, $\left\{ {{\lambda _i}} \right\}_{i = 1}^M$ and $\left\{ {{\omega _i}\left( x \right)} \right\}_{i = 1}^M$ are eigenvalues and eigenfunctions of the covariance function ${C_{\omega \omega }}\left( {{x_1},{x_2}} \right)$ of the random field $\omega \left( {x,\theta } \right)$, and they are solutions of the homogenous Fredholm integral equation of the second kind~\citep{ghanem2003stochastic, phoon2002simulation},
\begin{equation} \label{KLIE}
  \int_D {{C_{\omega \omega }}\left( {{x_1},{x_2}} \right){\omega _i}\left( {{x_1}} \right)d{x_1}}  = {\lambda _i}{\omega _i}\left( {{x_2}} \right)
\end{equation}
Due to the symmetry and the positive definiteness of covariance kernel ${C_{\omega \omega }}\left( {{x_1},{x_2}} \right)$, the eigenfunctions $\left\{ {{\omega _i}\left( x \right)} \right\}_{i = 1}^M$ form a complete orthogonal set satisfying the equation
\begin{equation} \label{Wi_O}
  \int_D {{\omega _i}\left( x \right){\omega _j}\left( x \right)dx}  = {\delta _{ij}}
\end{equation}
where $\delta _{ij}$ is the Kronecker delta function. An explicit expression for the random variables $\left\{ {{\xi _i}\left( \theta  \right)} \right\}_{i = 1}^M$ in Equation \eqref{RF} can be obtained by
\begin{equation} \label{KeSi}
  {\xi _i}\left( \theta  \right) = \frac{1}{{\sqrt {{\lambda _i}} }}\int_D {\left[ {\omega \left( {x,\theta } \right) - {\omega _0}\left( x \right)} \right]{\omega _i}\left( x \right)dx}
\end{equation}
which is a set of uncorrelated standardized random variables and satisfy
\begin{equation} \label{KeSi_O}
  E\left\{ {{\xi _i}\left( \theta  \right)} \right\} = 0,~E\left\{ {{\xi _i}\left( \theta  \right){\xi _j}\left( \theta  \right)} \right\} = {\delta _{ij}}
\end{equation}
where $E\left\{  \cdot  \right\}$ is the expectation operator.

The Karhunen-Lo\`{e}ve expansion \eqref{KL} offers a unified and powerful tool~\citep{phoon2002simulation, phoon2005simulation} for representing stationary and nonstationary, Gaussian and non-Gaussian random fields with explicitly known covariance functions. Karhunen-Lo\`{e}ve expansion is optimal among series expansion methods in the global mean square error with respect to the number of random variables in the representation, which means that only a few terms $M$ are required in order to capture most of randomness, thus it has received much attentions in many disciplines. In stochastic finite element analysis, it has been widely used to discretize the random fields representing the randomness of structures and excitations. It is worth mentioning that the implementation of Karhunen-Lo\`{e}ve expansion requires solutions of the integral equation \eqref{KLIE} with the covariance function as the integral kernel. Although only a limited number of analytical eigen-solutions are available~\citep{ghanem2003stochastic}, the solution of the integral equation can be numerically approximated for random fields with arbitrary covariance functions. For random fields that are defined on two- and three-dimensional domains, the finite element method becomes the only available method for the discretization of the multi-dimensional integral eigenvalue problems~\citep{2019An, zheng2017simulation}. In this paper, the generation of the finite element mesh for random fields is same to that for responses.

\subsection{Stochastic finite element equations} \label{sec22}
We simply recall the deterministic finite element method of the relevant formulation before dealing with stochastic problems. The deterministic finite element method in linear elasticity defined on $\Omega$ eventually derive a $N \times N$ linear system
\begin{equation} \label{FEM}
  Ku = F
\end{equation}
where $N$ is the number of degrees of freedom, $K$, $u$, $F$ are global stiffness matrix, displacement vector and load vector, respectively. By assembling the element stiffness matrices $k^e$, the global stiffness matrix $K$ can be obtained as
\begin{equation} \label{EleK}
  {k^e} = \int_{{\Omega _e}} {{B^T}DBd{\Omega _e}}
\end{equation}
where $B$ and $D$ stand for the strain matrix and the elasticity matrix, respectively.

We suppose that the material Young's modulus is a random field~\citep{sudret2000stochastic} and can be written as the form in Equation \eqref{KL},
\begin{equation} \label{DM}
  D\left( {x,\theta } \right) = {D_0}\left[ {{\omega _0}\left( x \right){\rm{ + }}\sum\limits_{i = 1}^M {{\xi _i}\left( \theta  \right)\sqrt {{\lambda _i}} {\omega _i}\left( x \right)} } \right]
\end{equation}
where $D_0$ is a constant matrix, and random variables $\left\{ {{\xi _i}\left( \theta  \right)} \right\}_{i = 1}^M$ construct a $M$-dimensional stochastic space. By substituting Equation \eqref{DM} into Equation \eqref{EleK}, the element stiffness matrix thus becomes as,
\begin{equation} \label{keR}
  {k^e}\left( \theta  \right) = k_0^e + \sum\limits_{i = 1}^M {{\xi _i}\left( \theta  \right)k_i^e}
\end{equation}
where $k_0^e$ is the mean element stiffness matrix given by
\begin{equation} \label{keR0}
  k_0^e = \int_{{\Omega _e}} {{\omega _0}\left( x \right){B^T}{D_0}Bd{\Omega _e}}
\end{equation}
and $k_i^e$ are deterministic matrices given by
\begin{equation} \label{keRi}
  k_i^e = \int_{{\Omega _e}} {\sqrt {{\lambda _i}} {\omega _i}\left( x \right){B^T}{D_0}Bd{\Omega _e}}
\end{equation}
The stochastic global stiffness matrix $K\left( \theta  \right)$ in the stochastic finite element equation \eqref{eq:SFEM} is obtained by assembling the stochastic element stiffness matrices ${k^e}\left( \theta  \right)$,
\begin{equation} \label{KGR}
  K\left( \theta  \right) = \sum\limits_{i = 0}^M {{\xi _i}\left( \theta  \right){K_i}}
\end{equation}
where ${\xi _0}\left( \theta  \right) \equiv 1$ and global matrices $K_i$ are obtained by assembling element matrices $k_i^e$ in the way similar to the deterministic case. In a similar way, we can get the stochastic global load vector as
\begin{equation} \label{FGR}
  F\left( \theta  \right) = \sum\limits_{l = 0}^Q {{\eta _l}\left( \theta  \right){F_l}}
\end{equation}
After assembling the stochastic global stiffness matrix $K\left( \theta  \right)$ and the stochastic global load vector $F\left( \theta  \right)$, the stochastic finite element equation \eqref{eq:SFEM} can be rewritten as
\begin{equation} \label{SFE0}
  \left( {\sum\limits_{i = 0}^M {{\xi _i}\left( \theta  \right){K_i}} } \right)u\left( \theta  \right) = \sum\limits_{l = 0}^Q {{\eta _l}\left( \theta  \right){F_l}}
\end{equation}

The high precision solution of Equation \eqref{SFE0} is one of the most important problems of the stochastic finite element method. Spectral stochastic finite element method (SSFEM)~\citep{xiu2010numerical, ghanem2003stochastic} is a popular method in the past few decades. SSFEM represents the stochastic response $u\left( \theta  \right)$ through polynomial chaos expansion (PCE) and transform Equation \eqref{SFE0} into a deterministic finite element equation by stochastic Galerkin projection. The size of the deterministic finite element equation depends directly on the number of terms retained in the PCE and the number of degrees of freedom $N$, and the computational cost for the solution of this system is much larger than that of the original problem.
Although several improved methods~\citep{Ghanem1996Precondition, nouy2010proper} have been developed to decrease computational costs, the Curse of Dimensionality still limit SSFEM to low-dimensional stochastic problems, thus it is crucial to develop a new method for the solution of Equation \eqref{SFE0}.

\section{A new method for solving stochastic finite element equations}\label{sec3}
In order to avoid the difficulties of SSFEM, in this section, we propose a new method for solving the sochastic finite element equation \eqref{SFE0} defined in low- and high-dimensional stochastic spaces. A natural idea is to represent the stochastic solution $u\left( \theta  \right)$ of Equation \eqref{SFE0} by use of random field expansions, however common methods are inactive since we almost know nothing about $u\left( \theta  \right)$ except the governing equation \eqref{SFE0}. Inspired by Karhunen-Lo\`{e}ve expansion \eqref{KL} and the general spectral decomposition~\citep{nouy2007generalized}, we construct the $u\left( \theta  \right)$ as
\begin{equation} \label{u0}
  u\left( \theta  \right) = \sum\limits_{i = 1}^\infty  {{\lambda _i}\left( \theta  \right){d_i}}
\end{equation}
where $\left\{ {{\lambda _i}(\theta )} \right\}_{i = 1}^\infty$ are random variables and $\left\{ {{d_i}} \right\}_{i = 1}^\infty$ are deterministic discretized basis vectors. Similar to the orthogonal conditions Equation \eqref{Wi_O} and \eqref{KeSi_O} of Karhunen-Lo\`{e}ve expansion, the following bi-orthogonal condition is introduced
\begin{equation} \label{OrC}
  d_i^T{d_j} = {\delta _{ij}},~E\left\{ {{\lambda _i}\left( \theta  \right){\lambda _j}\left( \theta  \right)} \right\} = {\kappa _i}{\delta _{ij}}
\end{equation}
where $E\{  \cdot \}$ is the expectation operator and ${\kappa _i} = E\{ {\lambda _i^2( \theta )} \}$.

It is shown in expansion \eqref{u0} that the solution space of $u\left( \theta  \right)$ is decoupled into a stochastic space and a deterministic space and it allows to compute $\left\{ {{\lambda _i}(\theta )} \right\}_{i = 1}^\infty$ in the stochastic space and  $\left\{ {{d_i}} \right\}_{i = 1}^\infty$ in the deterministic space, respectively. In this way, the difficulties in expanding the unknown solution random field of Equation \eqref{eq:SFEM} can be overcome. One only requires to seek a set of deterministic orthogonal vectors $\left\{ {{d_i}} \right\}_{i = 1}^\infty$ and the corresponding uncorrelated random variables $\left\{ {{\lambda _i}(\theta )} \right\}_{i = 1}^\infty$ such that the expanded solution in Equation \eqref{u0} satisfies the Equation \eqref{eq:SFEM}. In practical, we truncate Equation \eqref{u0} at the $k$-th term as,
\begin{equation} \label{uk}
  {u_k}\left( \theta  \right) = \sum\limits_{j = 1}^k {{\lambda _j}\left( \theta  \right){d_j}}
\end{equation}
As mentioned above, neither $\left\{ {{d_i}} \right\}_{i = 1}^k$ nor $\left\{ {{\lambda _i}(\theta )} \right\}_{i = 1}^k$ is known a priori, a natural choice is to successively determine these unknown couples $\left\{ {{\lambda _i}\left( \theta  \right), {d_i}} \right\}$ one after another via iterative methods. In order to compute the couple $\left( {{\lambda _k}\left( \theta  \right),{d_k}} \right)$, we suppose that the approximate solution $u_{k-1}(\theta)$ has been obtained, then substituting Equation \eqref{uk} into Equation \eqref{eq:SFEM} yields,
\begin{equation} \label{Kuk}
  K\left( \theta  \right)\left[ {\sum\limits_{j = 1}^{k - 1} {{\lambda _j}\left( \theta  \right){d_j}}  + {\lambda _k}\left( \theta  \right){d_k}} \right] = F\left( \theta  \right)
\end{equation}

If random variable ${\lambda _k}(\theta )$ has been determined (or given an initial value), the $d_k$ can be determined using stochastic Galerkin method and a dedicated iteration~\citep{nouy2007generalized}, this corresponds
\begin{equation} \label{Elambda}
  E\left\{ {{\lambda _k}(\theta )K(\theta )\left[ {\sum\limits_{j = 1}^{k - 1} {{\lambda _j}\left( \theta  \right){d_j}}  + {\lambda _k}(\theta ){d_k}} \right]} \right\} = E\left\{ {{\lambda _k}(\theta )F(\theta )} \right\}
\end{equation}
Considering Equation \eqref{KGR} and \eqref{FGR}, the Equation \eqref{Elambda} about $d_k$ can be simplified as,
\begin{equation} \label{Eq_d}
  \left( {\sum\limits_{i = 0}^M {{c_{ikk}}{K_i}} } \right){d_k} = \sum\limits_{l = 0}^Q {{b_{kl}}{F_l}}  - \sum\limits_{i = 0}^M {\sum\limits_{j = 1}^{k - 1} {{c_{ijk}}{K_i}{d_j}} }
\end{equation}
where
\begin{equation} \label{Cijk}
  {c_{ijk}} = E\left\{ {{\xi _i}\left( \theta  \right){\lambda _j}\left( \theta  \right){\lambda _k}\left( \theta  \right)} \right\}, \; {b_{kl}} = E\left\{ {{\eta _l}\left( \theta  \right){\lambda _k}\left( \theta  \right)} \right\}
\end{equation}

Once $d_k$ has been determined in Equation \eqref{Eq_d}, the random variable ${\lambda _k}(\theta )$ can be subsequently updated via the similar procedure. This requires to multiply ${d_k}$ on both sides of Equation \eqref{Kuk} to yield
\begin{equation} \label{Mul_dk}
  d_k^TK\left( \theta  \right)\left[ {\sum\limits_{j = 1}^{k - 1} {{\lambda _j}\left( \theta  \right){d_j}}  + {\lambda _k}\left( \theta  \right){d_k}} \right] = d_k^TF\left( \theta  \right)
\end{equation}
Considering Equation \eqref{KGR} and \eqref{FGR}, the Equation \eqref{Mul_dk} about ${\lambda _k}(\theta )$ can be simplified as,
\begin{equation} \label{Eq_lambda}
  \left( {\sum\limits_{i = 0}^M {{g_{ikk}}{\xi _i}\left( \theta  \right)} } \right){\lambda _k}\left( \theta  \right) = \sum\limits_{l = 0}^Q {{h_{kl}}{\eta _l}\left( \theta  \right)}  - \sum\limits_{i = 0}^M {\sum\limits_{j = 1}^{k - 1} {{g_{ijk}}{\xi _i}\left( \theta  \right){\lambda _j}\left( \theta  \right)} }
\end{equation}
where
\begin{equation} \label{Gijk}
  {g_{ijk}} = d_k^T{K_i}{d_j}, \; {h_{kl}} = d_k^T{F_l}
\end{equation}

The classical SSFEM is to represent the stochastic solution of nodes $\left\{ {{u_i}\left( \theta  \right)} \right\}_{i = 1}^N$ in terms of a set of polynomial chaos and transforms the original stochastic finite element equation into a deterministic finite element equation with size $N \times \frac{{\left( {M + p} \right)!}}{{M!p!}}$, where $\left(  \cdot  \right)!$ represents the factorial operator, $N$, $M$ and $p$ are the number of system degrees of freedom, the number of random variables and the order of polynomial chaos expansion, respectively. The size of the deterministic finite element equation is significantly higher than that of the original stochastic finite element equation. For instance, the size is $1 \times 10^6$ when $N = 1000$, $M = 10$ and $p = 4$, which leads to the Curse of Dimensionality, and is prohibitive for problems with high stochastic dimensions and large scales.
		
The method in this paper decouples the original stochastic finite element equation into a deterministic finite element equation \eqref{Eq_d} with size $N$ and one-dimensional stochastic algebraic equation \eqref{Eq_lambda}. The iteration process of Equation \eqref{Elambda} and \eqref{Mul_dk} was used in the paper~\citep{nouy2007generalized} to solve linear stochastic partial differential equations and subsequently applied to time-dependent and nonlinear problems~\citep{nouy2009generalized, nouy2010priori}. The key of this method is to transform the original SPDE to a deterministic PDE and a stochastic algebraic equation like Equation \eqref{Eq_lambda}. The method for solving Equation \eqref{Eq_lambda}-like is to represent the random variable $\lambda_k(\theta)$ in terms of a set of polynomial chaos and transforms Equation \eqref{Eq_lambda}-like into a deterministic equation with size $\frac{{\left( {M + p} \right)!}}{{M!p!}}$, which greatly alleviates the Curse of Dimensionality, but is still prohibitive for problems with high stochastic dimensions.
In order to avoid this difficulty, we develop a simulation method to determine $\lambda_k(\theta)$. For each realization of $\{ {\theta ^{(r)}}\} _{r = 1}^R$, the ${\lambda _k}({\theta ^{(r)}})$ can be obtained by solving \eqref{Eq_lambda} as,
\begin{equation} \label{Sam_lambda}
  {\lambda _k}\left( {{\theta ^{\left( r \right)}}} \right) = \frac{{\sum\limits_{l = 0}^Q {{h_{kl}}{\eta _l}\left( {{\theta ^{\left( r \right)}}} \right)}  - \sum\limits_{i = 0}^M {\sum\limits_{j = 1}^{k - 1} {{g_{ijk}}{\xi _i}\left( {{\theta ^{\left( r \right)}}} \right){\lambda _j}\left( {{\theta ^{\left( r \right)}}} \right)} } }}{{\sum\limits_{i = 0}^M {{g_{ikk}}{\xi _i}\left( {{\theta ^{\left( r \right)}}} \right)} }}
\end{equation}

It is important to note that Equation \eqref{Sam_lambda} has become a one-dimensional linear algebraic equation about ${\lambda _k}({\theta ^{(r)}})$. Compared to classic methods, we do not need to choose the type and order of polynomial chaos. The total computational cost for determining $\{ {\lambda _k}({\theta ^{(r)}})\} _{r = 1}^R$ is very low even for high stochastic dimensions, which hopefully avoid the Curse of Dimensionality. Then statistical methods are readily introduced to obtain $\lambda_k(\theta)$ from samples $\{ {\lambda _k}({\theta ^{(r)}})\} _{r = 1}^R$. Hence, this method will be particularly appropriate for a wide class of high-dimensional stochastic problems in practice.

\section{Algorithm implementation}\label{sec4}

\begin{algorithm}[h]
	\caption{Algorithm for solving linear stochastic finite element equations} \label{alg1}
	\begin{algorithmic}[1]
		\While {${\varepsilon _{global}} > {\varepsilon _1}$}
		\label{step1}
		\State initialize $\lambda _k^{(0)}(\theta )$
		\label{step2}
		\While {${\varepsilon _{local}} > {\varepsilon _2}$}
		\label{step3}
		\State compute $d_k^{\left( j \right)}$ by solving Equation \eqref{Eq_d}
		\label{step4}
		\State orthogonalization $d_k^{\left( j \right)} \bot {d_i},\; i = 1, \cdots, k - 1$ by Equation \eqref{OrDL} and normalization $d_k^{\left( j \right)} = \frac{{d_k^{\left( j \right)}}}{{\left\| {d_k^{\left( j \right)}} \right\|}}$
		\label{step5}
		\State compute $\lambda _k^{(j)}(\theta )$ by Equation \eqref{Sam_lambda}
		\label{step6}
		\State orthogonalization $\lambda _k^{\left( j \right)}\left( \theta  \right) \bot {\lambda _i}\left( \theta  \right),\; i = 1, \cdots, k - 1$ by Equation \eqref{OrDL} 
		\label{step7}
		\State compute local error ${\varepsilon _{local}}$, $j=j+1$
		\label{step8}
		\EndWhile
		\label{step9}
		\State update $u(\theta )$ as ${u_k}(\theta ) = \sum\limits_{i = 1}^{k - 1} {{\lambda _i}(\theta ){d_i}}  + {\lambda _k}(\theta ){d_k}$
		\label{step10}
		\State compute global error ${\varepsilon _{global}}$, $k=k+1$
		\label{step11}
		\EndWhile
		\label{step12}
	\end{algorithmic}
\end{algorithm}

The above procedure for solving the stochastic finite element equation \eqref{eq:SFEM} is summarized in Algorithm \ref{alg1}, which consists of a outer loop procedure and a inner loop procedure. The inner loop, which is from step \ref{step3} to \ref{step9}, is used to determine the couple of $(\lambda_k(\theta), d_k)$. With an initial random variable $\lambda _k^{\left( 0 \right)}\left( \theta  \right)$ given in step \ref{step2}, $d_k^{\left( j \right)}$ can be determined in step \ref{step4} and \ref{step5}, where superscript $j$ represents the \emph{j}-th round of iteration. With the obtained $d_k^{\left( j \right)}$, the random variable $\lambda _k^{\left( j \right)}\left( \theta  \right)$ is then updated in step \ref{step6} and \ref{step7}. While the outer loop, which is from step \ref{step1} to \ref{step12}, corresponds to recursively building the set of couples and  then generates a set of couples such that the approximate solution in step \ref{step10} satisfies Equation \eqref{eq:SFEM}.

Note that both $d_k^{\left( j \right)}$ and $\lambda _k^{\left( j \right)}\left( \theta  \right)$ require orthogonalizations such that the bi-orthogonal condition in Equation \eqref{OrC} holds along the whole process, here we use the Gram-Schmidt Orthogonalization method in step \ref{step5} and \ref{step7}. It is written as,
\begin{equation} \label{OrDL}
\left\{ {\begin{array}{*{20}{l}}
	{d_k^{\left( j \right)} = d_k^{\left( j \right)} - \sum\limits_{i = 1}^{k - 1} {\left( {d_k^{\left( j \right)T}{d_i}} \right){d_i}} }\\
	{\lambda _k^{\left( j \right)}\left( \theta  \right) = \lambda _k^{\left( j \right)}\left( \theta  \right) - \sum\limits_{i = 1}^{k - 1} {\frac{{E\left\{ {\lambda _k^{\left( j \right)}\left( \theta  \right){\lambda _i}\left( \theta  \right)} \right\}}}{{E\left\{ {\lambda _i^2\left( \theta  \right)} \right\}}}{\lambda _i}\left( \theta  \right)} }
	\end{array}} \right.,k \ge 2
\end{equation}

For practical purposes, a certain number of truncated items are retained of the solution $u\left( \theta  \right)$. The truncation criterion in step \ref{step1} is considered as a global error, which is defined as,
\begin{eqnarray} \label{err_glo}
  {\varepsilon _{global}} &= &\frac{{E\left\{ {\Delta u_k^2\left( \theta  \right)} \right\}}}{{E\left\{ {u_k^2\left( \theta  \right)} \right\}}} = \frac{{E\left\{ {\lambda _k^2\left( \theta  \right)} \right\}d_k^T{d_k}}}{{\sum\limits_{i = 1}^k {\sum\limits_{j = 1}^k {E\left\{ {{\lambda _i}\left( \theta  \right){\lambda _j}\left( \theta  \right)} \right\}d_i^T{d_j}} } }}  \nonumber\\
  &= &\frac{{E\left\{ {\lambda _k^2\left( \theta  \right)} \right\}}}{{\sum\limits_{i = 1}^k {E\left\{ {\lambda _i^2\left( \theta  \right)} \right\}} }}
\end{eqnarray}
which measures the contribution of the $k$-th couple $\left( {{\lambda _k}\left( \theta  \right),{d_k}} \right)$ to the stochastic solution $u\left( \theta  \right)$ and converges to the final solution when it achieves the required precision. Further, the stop criterion for computing each couple $\left( {{\lambda _k}\left( \theta  \right),{d_k}} \right)$ is considered as a local error and defined as,
\begin{eqnarray} \label{err_loc}
  {\varepsilon _{local}} = \frac{{\left\| {d_k^{(j)} - d_k^{(j - 1)}} \right\|}}{{\left\| {d_k^{(j - 1)}} \right\|}} = \left\| {d_k^{(j)} - d_k^{(j - 1)}} \right\|
\end{eqnarray}
which measures the difference between $d_k^{(j)}$ and $d_k^{(j-1)}$ and the calculation is stopped when $d_k^{(j)}$ is almost the same as $d_k^{(j-1)}$.

\section{Applications}\label{sec5}
The numerical implementation of the proposed method is illustrated with the aid of three practical applications. The first application consists of an electric pylon frame with stochastic material properties and a stochastic load. The second application is a roof truss under stochastic wind loads defined in low-dimensional and high-dimensional stochastic spaces, which is to illustrate the efficiency of applying the proposed method to high-dimensional stochastic problems. The third example tests the ability of the proposed method for dealing with a large-scale engineering problem given by computing the deformation of a tunnel under the action of self-weight. These three exampls serve to verify the validity and accuracy of the proposed method and demonstrate that there is the same solution construct for different problems.


\subsection{Response of electric pylon frame with stochastic material property}\label{Example-1}

\begin{figure}
	\centerline{\includegraphics[width=0.4\linewidth]{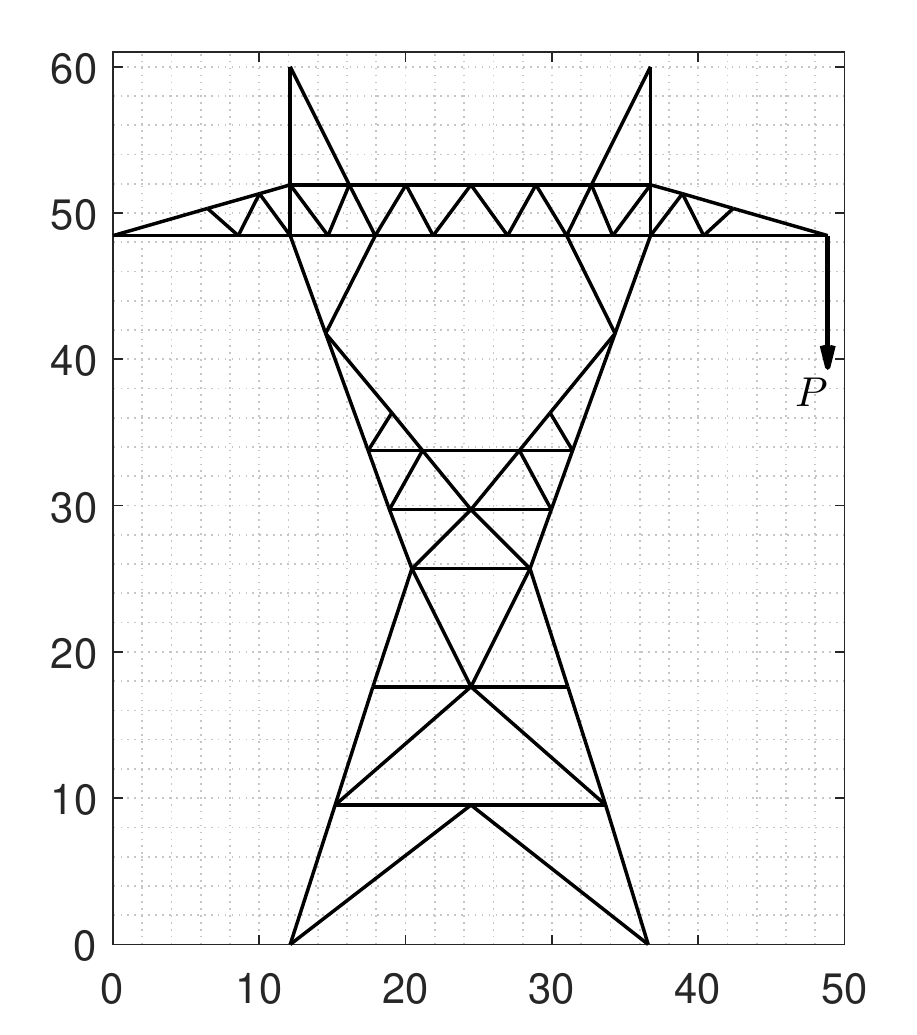}}
	\caption{91-element electric pylon frame \label{fig_e1_01}}
\end{figure}

In this problem, we consider a frame system as shown in Figure \ref{fig_e1_01}, which is a electric pylon frame consisting of 91 elements with square cross-sections. A load $P$ is applied vertically downward at the far right tip of the arm of the pylon. Clamped boundary conditions are applied at the base of the frame model. Spatial nodes of the electric pylon frame model are defined in Table \ref{tab_e1}. All elements of the electric pylon frame are constructed of 300M steel and have identical cross-sectional areas. Deterministic material properties are given as, mass density $\rho = 7.8\rm{{g \mathord{\left/{\vphantom {g {c{m^3}}}} \right.\kern-\nulldelimiterspace} {c{m^3}}}}$, cross-sectional area $\overline{A}=4\rm{c{m^2}}$, Young's modulus $\overline{E} = 200\rm{GPa}$.

\begin{table}
	\caption{Nodal definitions of the electric pylon frame \label{tab_e1}}
	\centering
	\begin{tabular}{rrrrrrrrrrrr}
		\toprule
		\bf{node} & $x$ & $y$ & \bf{node} & $x$ & $y$ \bf{node} & $x$ & $y$ & \bf{node} & $x$ & $y$ \\
		\midrule
		\bf{1}   & 12.11 &	0.00   &  \bf{13} & 29.96 &  29.73 & \bf{25} &	14.70 &	48.45  & \bf{37} &	38.92 &	51.32\\
		\bf{2}   & 36.58 &	0.00   & \bf{14} & 17.44 &  33.76 &  \bf{26} &	17.93 &	48.45 &  \bf{38} &	12.11 &	51.92 \\
		\bf{3}   &	 15.18 &	9.53   &  \bf{15} & 21.16 &  33.76 &  \bf{27} &	21.88 &	48.45 & \bf{39} &	16.15 &	51.92\\
		\bf{4}   &	 24.47 &	9.53   &  \bf{16} &	 27.78 &	33.76 & \bf{28} &	26.97 &	48.45 & \bf{40} &	20.03 &	51.92 \\
		\bf{5}   &	 33.67 &	9.53   &  \bf{17} &	 31.41 &	33.76 & \bf{29} &	31.01 &	48.45 &\bf{41} &	24.47 &	51.92 \\ 
		\bf{6}   & 17.77 &	17.61 & \bf{18} &	 19.06 &	36.34 &   \bf{30} &	34.16 &	48.45 & \bf{42} &	28.91 &	51.92\\
		\bf{7}   & 24.47 &	17.61 &  \bf{19} &	 29.88 &	36.34 & \bf{31} &	36.74 &	48.45 & \bf{43} &	32.71 &	51.92 \\
		\bf{8}   &	 31.09 &	17.61 & \bf{20} &	 14.54 &	41.77 & \bf{32} &	40.38 &	48.45 & \bf{44} &	36.74 &	51.92   \\
		\bf{9}   &	 20.43 &	25.69 & \bf{21} &	 34.32 &	41.77 & \bf{33} &	48.86 &	48.45 & \bf{45} &	12.11 &	60.00  \\
		\bf{10} & 28.51 &  25.69 &  \bf{22} &	 0.00   &	48.45 &  \bf{34} &	6.46   &	50.31 & \bf{46}	&   36.74 &	60.00\\ 
		\bf{11} & 18.90 &  29.73 & \bf{23} &	 8.56   &	48.45 &   \bf{35} &	42.40 &	50.31 & \\
		\bf{12} & 24.47 &  29.73 &\bf{24} &	12.11 &	48.45 & \bf{36} &	10.01 &	51.32 &  \\
		\bottomrule
	\end{tabular}
\end{table}

The response of the electric pylon forced under a load $P$ deeply depends on these parameters. In order to better reflect the structural response influenced by material and load variabilities, we consider the stochastic tensile stiffness and the stochastic bending stiffness as,
\begin{equation} \label{EA_EI}
  EA = \left( {{\xi _1}\left( \theta  \right) + 0.2{\xi _2}\left( \theta  \right)} \right)\overline {EA},~ 
  EI = \left( {{\xi _3}\left( \theta  \right) + 0.2{\xi _4}\left( \theta  \right)} \right)\overline {EI}
\end{equation}
and consider a stochastic load $P\left( \theta  \right)$ as
\begin{equation} \label{StoP}
  P\left( \theta  \right) = \left( {1 + {\xi _5}\left( \theta  \right) + {\xi _6}\left( \theta  \right)} \right)\overline{P}
\end{equation}
where $\overline{P} = 1000\rm{N}$. Indepdent random variables $\left\{ {{\xi _i}\left( \theta  \right)} \right\}_{i = 1}^6$ in Equation \eqref{EA_EI} and \eqref{StoP} satisfy
\begin{equation} \label{RVs}
  \log \left\{ {{\xi _i}\left( \theta  \right)} \right\}_{i = 1}^4 \sim N\left( {0,0.3} \right),~{\xi _5}\left( \theta  \right),{\xi _6}\left( \theta  \right) \sim N\left( {0,0.1} \right)
\end{equation}

Similar to the derivation of Equation \eqref{SFE0}, a stochastic finite element equation for this problem can be obtained as,
\begin{equation} \label{E1_SFEM}
\left( {\sum\limits_{i = 1}^4 {{\xi _i}\left( \theta  \right){K_i}} } \right)u\left( \theta  \right) = \left( {1 + \sum\limits_{i = 5}^6 {{\xi _i}\left( \theta  \right)} } \right)F
\end{equation}

\begin{figure}
	\centerline{\includegraphics[width=0.6\linewidth]{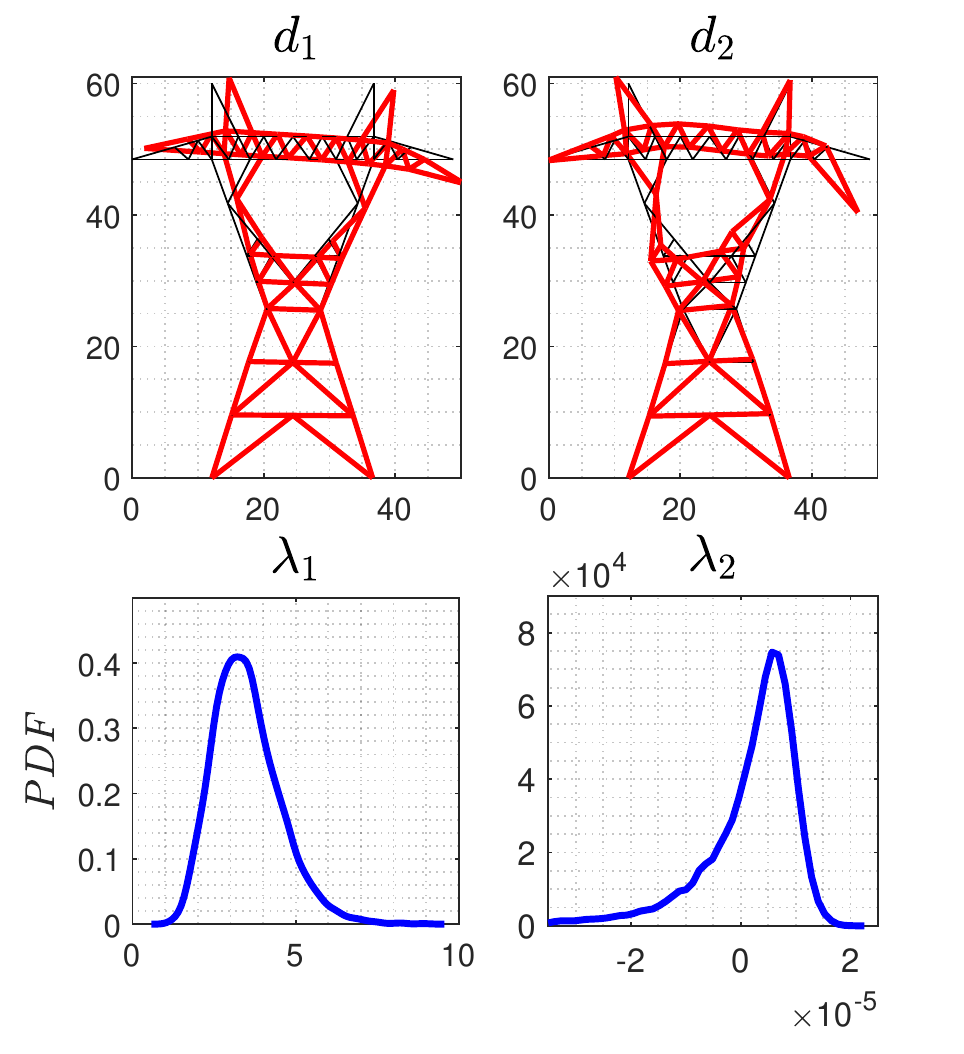}}
	\caption{Displacement components $\left\{ {{d_i}} \right\}_{i = 1}^2$ (top) and PDFs of corresponding random variables $\left\{ {{\lambda _i}\left( \theta  \right)} \right\}_{i = 1}^2$ (bottom) \label{E1_DL}}
\end{figure}

\begin{figure}
	\centerline{\includegraphics[width=0.6\linewidth]{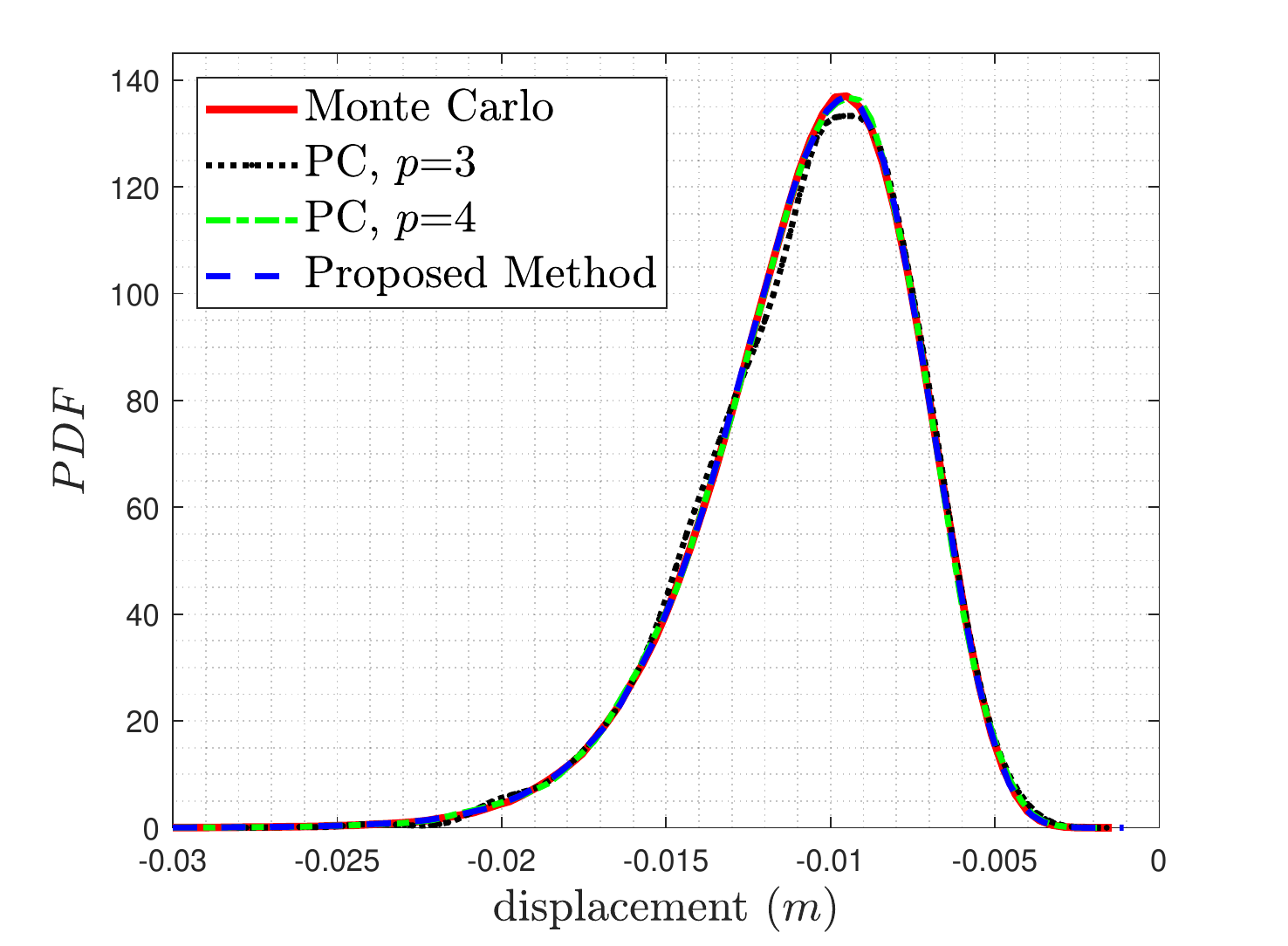}}
	\caption{Comparison of PDFs between the Monte Carlo simulation and the proposed method \label{E1_pdf}}
\end{figure}

In order to solve Equation \eqref{E1_SFEM} by use of Algorithm \ref{alg1}, the convergence criterias are set as $\varepsilon _1 = \varepsilon _2 = 10^{ - 6}$, $R = 1 \times 10^4$ random samples, i.e. $\left\{ {{\xi _i}\left( {{\theta ^{\left( r \right)}}} \right)} \right\}_{r = 1}^{1 \times {{10}^4}},~i = 1, \cdots ,2$, are given in Equation \eqref{E1_SFEM} and $\{ \lambda _k^{(0)}({\theta ^{(r)}})\} _{r = 1}^{1 \times {{10}^4}}$ are adopt in step \ref{step2} in Algorithm \ref{alg1}. In this example, only two retained terms in Equation \eqref{uk}, the displacement components $\left\{ {{d_i}\left( {x,y} \right)} \right\}_{i = 1}^2$ and probability density functions (PDFs) of corresponding random variables $\left\{ {{\lambda _i}\left( \theta  \right)} \right\}_{i = 1}^2$ shown in Figure \ref{E1_DL}, can achieve the required precision, which demonstrates the high efficiency of the proposed method. It is seen from Figure \ref{E1_DL} that, the second random variable ${\lambda _2}\left( \theta  \right)$ is very small and makes almost no contributions to the approximate solution ${u_2}\left( \theta  \right)$. In practical, we determine whether step \ref{step2} in Algorithm \ref{alg1} converges through computing the second term, thus the displacement component $d_2$ and the random variable ${\lambda _2}\left( \theta  \right)$ are necessary.

Further we compare the proposed method with existing methods, including Monte Carlo simulation~\citep{Papadrakakis1996Robust} and spectral stochastic finite element method (SSFEM)~\citep{xiu2010numerical, ghanem2003stochastic}. Here Hermite Polynomial Chaos (PC) of 6 standard Gaussian random variables are adopted in the SSFEM, and the order of PC is set as $p=3$ and $p=4$. We test the computational efficiency of these methods by use of a personal laptop (dual-core, Intel i7, 2.40GHz) and the computational times of the proposed method, PC ($p=3$), PC ($p=4$) and $1 \times {10^6}$ standard Monte Carlo simulations are 3.4s, 71.9s, 474.9s and 1412.7s, respectively, which demonstrates the high efficiency of the proposed method. 
Based on above methods, the resulted approximate PDFs of the response of the far right tip of the arm of the electric pylon are seen from Figure \ref{E1_pdf}. The result of the two-term approximation of the proposed method is in very good accordance with that from the Monte Carlo simulation, while the PC method requires fourth order ($p=4$) to achieve a similar accuracy. In addition, our method is based on random samples, thus can avoid choosing the order $p$ of PC basis. We observe in practice that the number of random samples has less influence on the computational cost. In the general case, the sample size is enough when it is sufficient to describe the statistical characteristics of the random variables.

\subsection{Response of roof truss under stochastic wind loads}\label{Example-2}
\begin{figure}[h]
	\centerline{\includegraphics[width=0.7\linewidth]{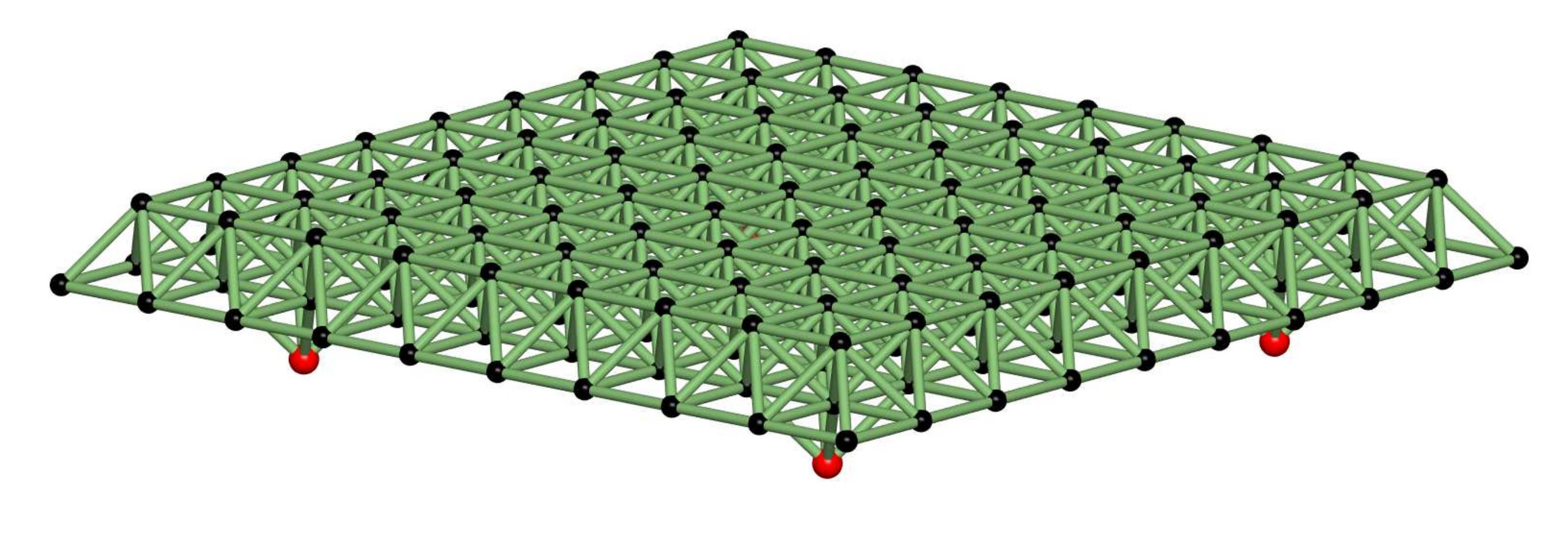}}
	\centerline{\includegraphics[width=0.7\linewidth]{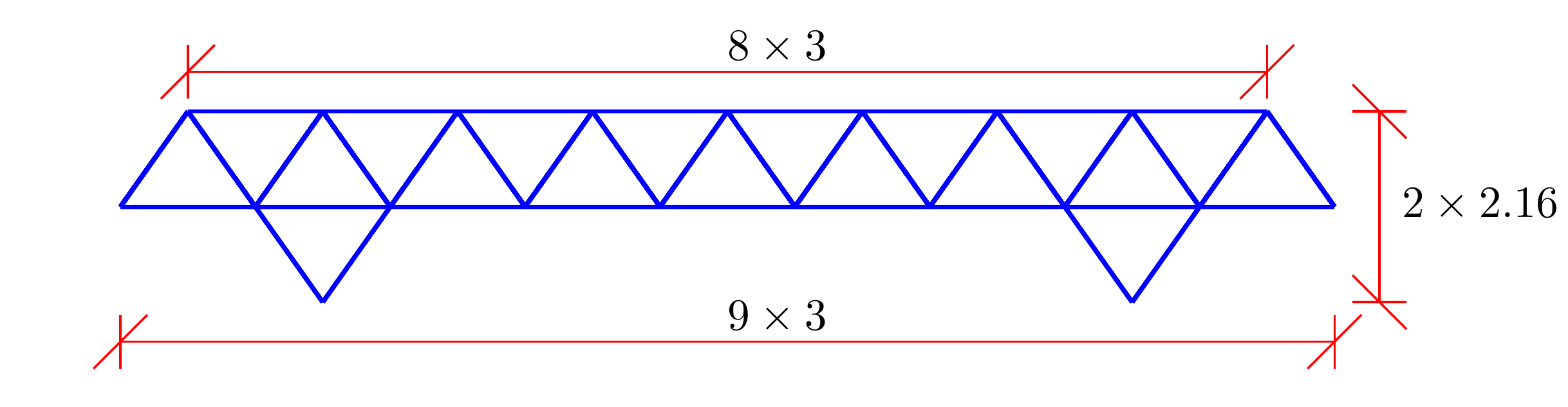}}
	\caption{Model of the roof truss\label{fig_e2_01}}
\end{figure}

In this example, we consider the stochastic response of a roof truss under a stochastic wind load acting vertically downward on the roof. The roof truss, as shown in Figure \ref{fig_e2_01}, includes 185 spatial nodes and 664 elements, where material properties of all members are set as Young's modulus $E = 209{\mathop{\rm GPa}\nolimits}$ and cross-sectional areas $A = 16{{\mathop{\rm cm}\nolimits} ^2}$.
The stochastic wind load is a random field with the covariance function ${C_{FF}}\left( {{x_1},{y_1};{x_2},{y_2}} \right) = \sigma _F^2{e^{ - {{\left| {{x_1} - {x_2}} \right|} \mathord{\left/ {\vphantom {{\left| {{x_1} - {x_2}} \right|} {{l_x}}}} \right. \kern-\nulldelimiterspace} {{l_x}}} - {{\left| {{y_1} - {y_2}} \right|} \mathord{\left/ {\vphantom {{\left| {{y_1} - {y_2}} \right|} {{l_y}}}} \right. \kern-\nulldelimiterspace} {{l_y}}}}}$, where the variance function $\sigma _F^2 = 0.15$, the correlation lengths ${l_x} = {l_y} = 24$, and it can be expanded by use of Karhunen-Lo\`{e}ve expansion Equation \eqref{KL} with a $M$-term truncated as
\begin{equation} \label{E2_RF}
  f\left( {x,y,\theta } \right) = \sum\limits_{i = 0}^{M} {{\xi _i}\left( \theta  \right){f_i}\left( {x,y} \right)}
\end{equation}
where ${\xi _0}\left( \theta  \right) \equiv 1$ and the mean function ${f_0}\left( {x,y} \right) = 10{\mathop{\rm kN}\nolimits}$. $\left\{ {{f_i}\left( {x,y} \right)} \right\}_{i = 1}^M$ is obtained by solving Equation \eqref{KL}. Based on the expansion Equation \eqref{E2_RF} of the stochastic wind load, the following stochastic finite element equation is obtained,
\begin{equation} \label{E2_SFEM}
  Ku\left( \theta  \right) = \sum\limits_{i = 0}^{M} {{\xi _i}\left( \theta  \right){F_i}}
\end{equation}

In this example, the initializations give the random samples $\left\{ {{\xi _i}\left( {{\theta ^{\left( r \right)}}} \right)} \right\}_{r = 1}^{1 \times {{10}^4}},~i = 1, \cdots ,M$ and the initial random variable samples $\{ \lambda _k^{(0)}({\theta ^{(r)}})\} _{r = 1}^{1 \times {{10}^4}}$, and set the convergence criterias as $\varepsilon _1 = \varepsilon _2 = 10^{ - 6}$. We first consider a low-dimensional case by choosing $M=10$. 
\begin{figure*}
	\centerline{\includegraphics[width=1.0\textwidth]{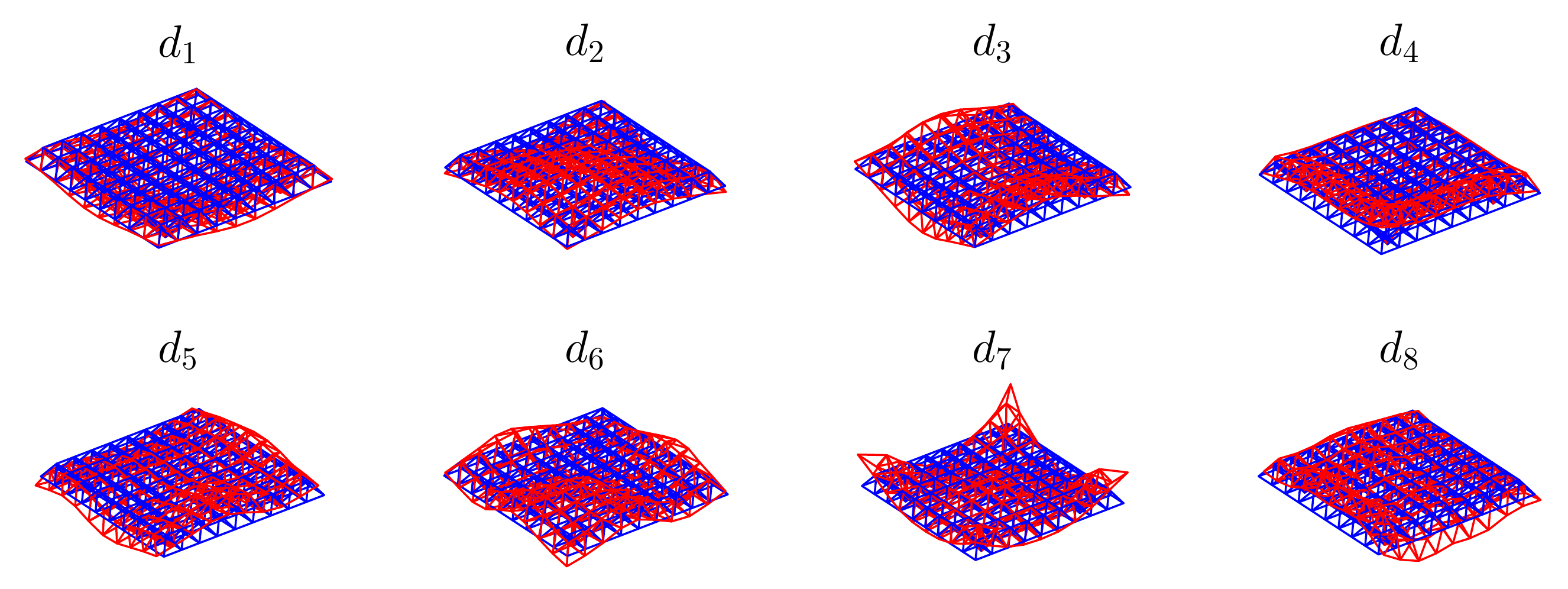}}
	\centerline{$\bf{(a)}$.~Displacement components $\left\{ {{d_i}} \right\}_{i = 1}^8$\label{fig_e2_021}}
	\begin{minipage}[t]{0.585\linewidth}
		\centerline{\includegraphics[width=1.0\linewidth]{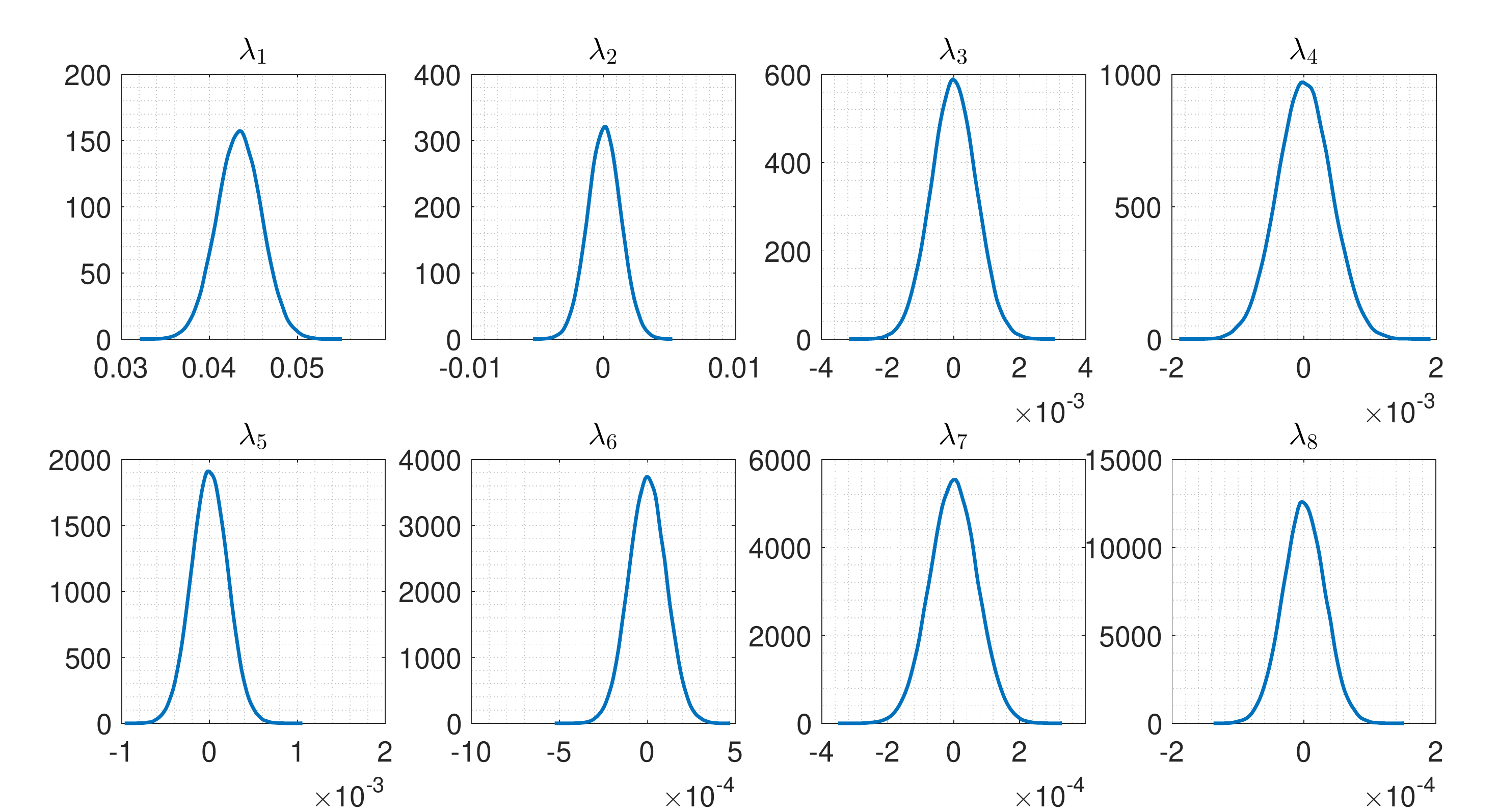}}
		\centerline{$\bf{(b)}$.~PDFs of $\left\{ {{\lambda _i}\left( \theta  \right)} \right\}_{i = 1}^8$ \label{fig_e2_022}}
	\end{minipage}
	\begin{minipage}[t]{0.415\linewidth}
		\centerline{\includegraphics[width=1.0\linewidth]{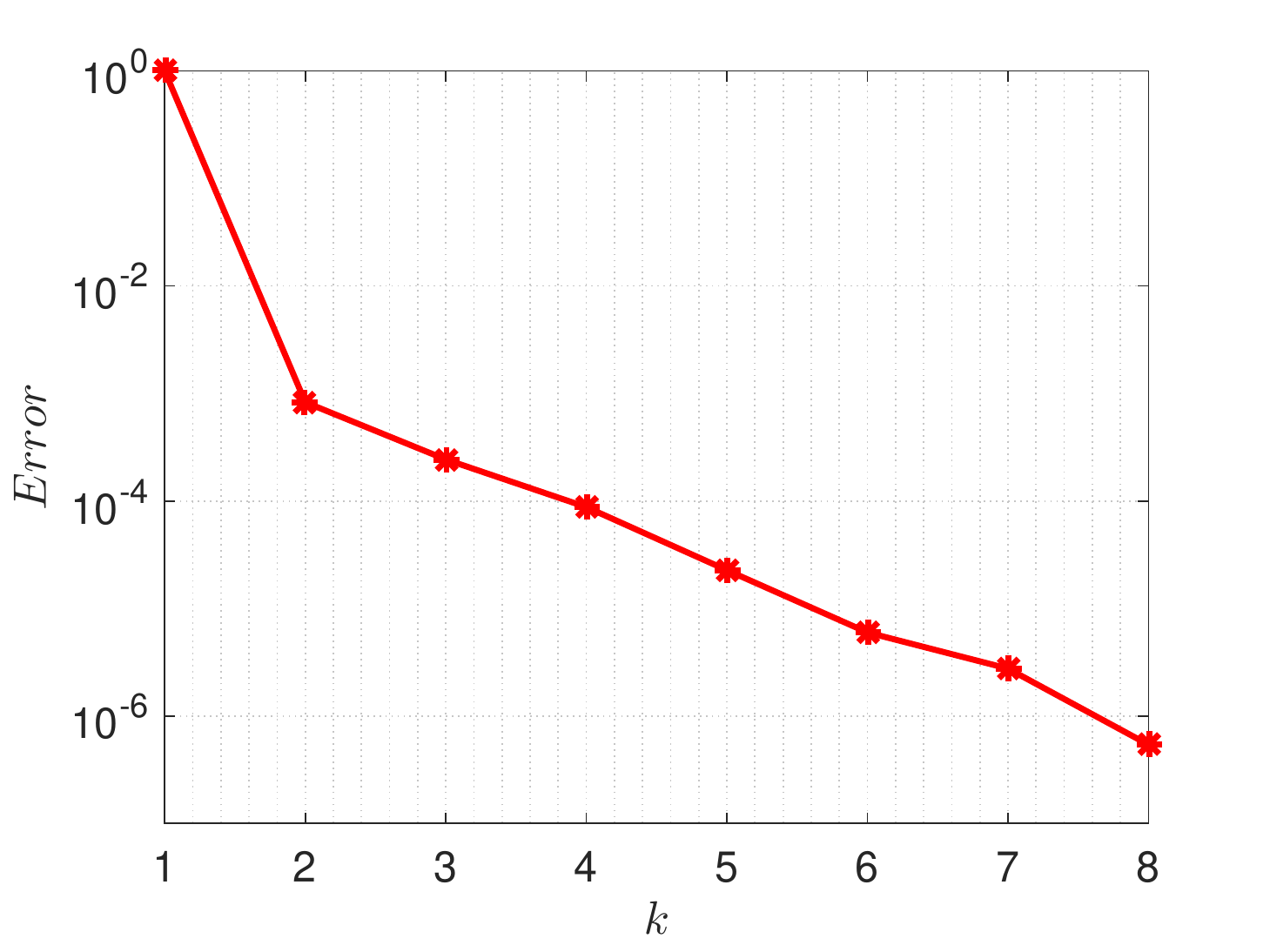}}
		\centerline{$\bf{(c)}$.~Iterative errors of $k$ retained items \label{fig_e2_023}}
	\end{minipage}
	\caption{Solutions of the couples $\left\{ {{\lambda _i}\left( \theta  \right),{d_i}} \right\}_{i = 1}^8$ and iterative errors of the solving process \label{fig_e2_02}}
\end{figure*}
It is seen from Figure \ref{fig_e2_02}c that the displacement components $\left\{ {{d_i}} \right\}$ and corresponding random variables $\left\{ {{\lambda _i}\left( \theta  \right)} \right\}$ can be determined after 8 iterations, which demonstrates the fast convergence rate of the proposed method. Correspondingly, the number of couples $\left( {{\lambda _k}\left( \theta  \right),{d_k}} \right)$ that constitute the stochastic response is adopted as $k=8$. As shown in Figure \ref{fig_e2_02}a and Figure \ref{fig_e2_02}b, with the increasing of the number of couples, the ranges of corresponding random variables are more closely approaching to zero, indicating that the contribution of the higher order random variables to the approximate solution decays dramatically. 

\begin{figure*}[hp]
	\begin{minipage}[t]{0.49\linewidth}
		\centerline{\includegraphics[width=1.0\linewidth]{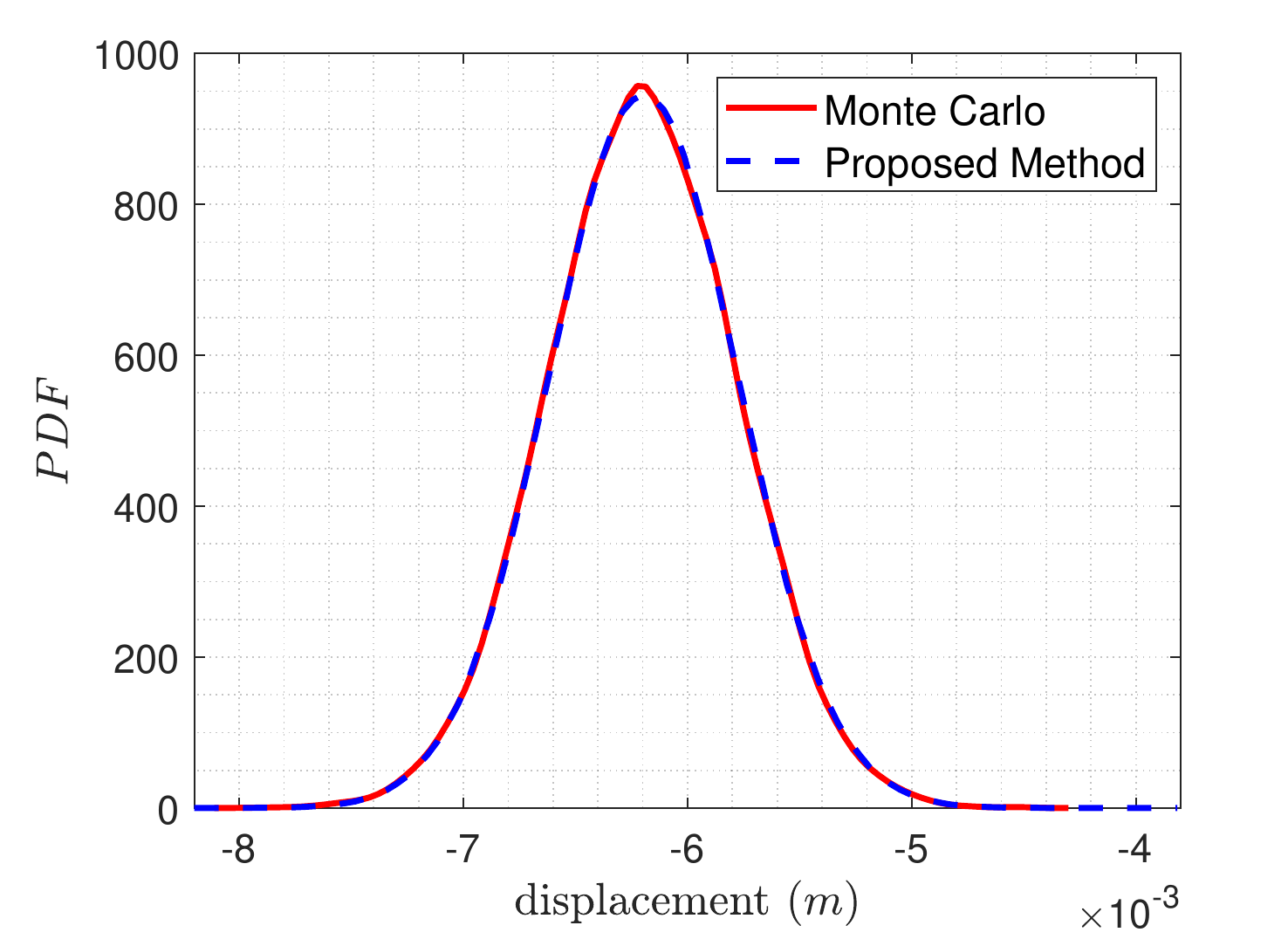}}
		\caption{Comparison of PDFs between the MCS and Algorithm \ref{alg1} \label{fig_e2_03}}
	\end{minipage}
	\begin{minipage}[t]{0.49\linewidth}
		\centerline{\includegraphics[width=1.0\linewidth]{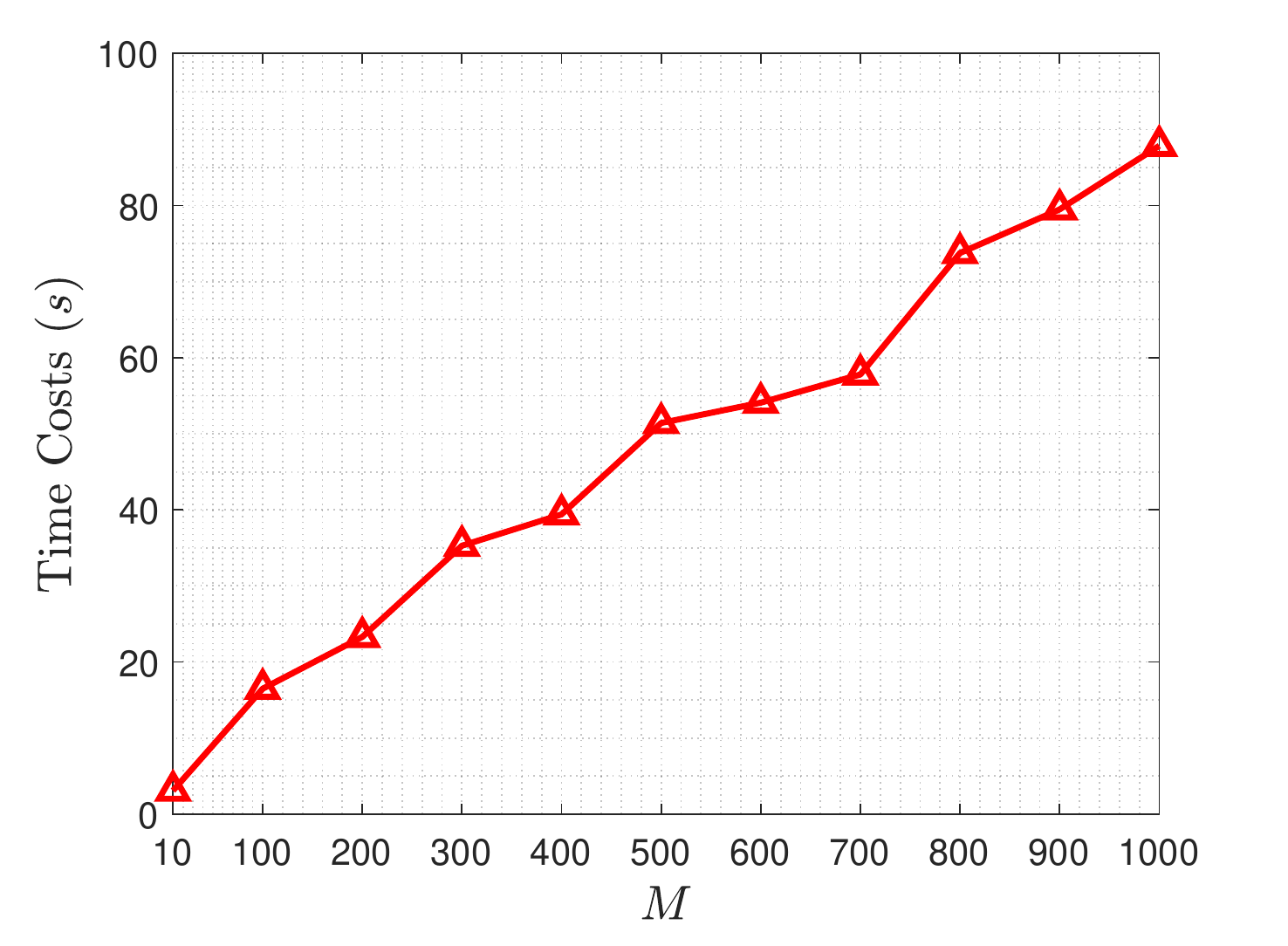}}
		\caption{Time costs of different stochastic dimensions $M = 10 \sim 1000$\label{E2_time}}
	\end{minipage}
\end{figure*}

For the maximum displacement of the whole roof truss, the resulted approximate PDF compared with $1 \times {10^6}$ standard Monte Carlo simulations (MCS) is seen in Figure \ref{fig_e2_03}, which indicates that the result of eight-term approximation is in very good accordance with that from the Monte Carlo simulation. According to our experience, further increasing the number of couples will not significantly improve the accuracy since the series in Equation \eqref{uk} has converged and thus the first few couples dominate the solution of the problem. This example demonstrates the success of our proposed construct of the stochastic solution and Algorithm \ref{alg1} for the solution of practical problems.

One of the main purposes the proposed method is to solve high-dimensional stochastic problems. Here we introduce the high-dimensional stochastic problems by choosing $M=100 \sim 1000$, and test the computational efficiency of different stochastic dimensions by use of a personal laptop (dual-core, Intel i7, 2.40GHz).
Computational costs for solving Equation \eqref{E2_SFEM} of different stochastic dimensions are shown in Figure \ref{E2_time}, which indicates that our proposed algorithm is efficient for high stochastic dimensions. The computational costs do not increase dramatically as the dimensions increase and is almost linear with the stochastic dimensions, which demonstrates the sucess of the proposed method for avoiding the Curse of Dimensionality.

\subsection{Deformation of tunnel under the action of self-weight}\label{Example-3}

\begin{figure}
	\centerline{\includegraphics[width=0.55\linewidth]{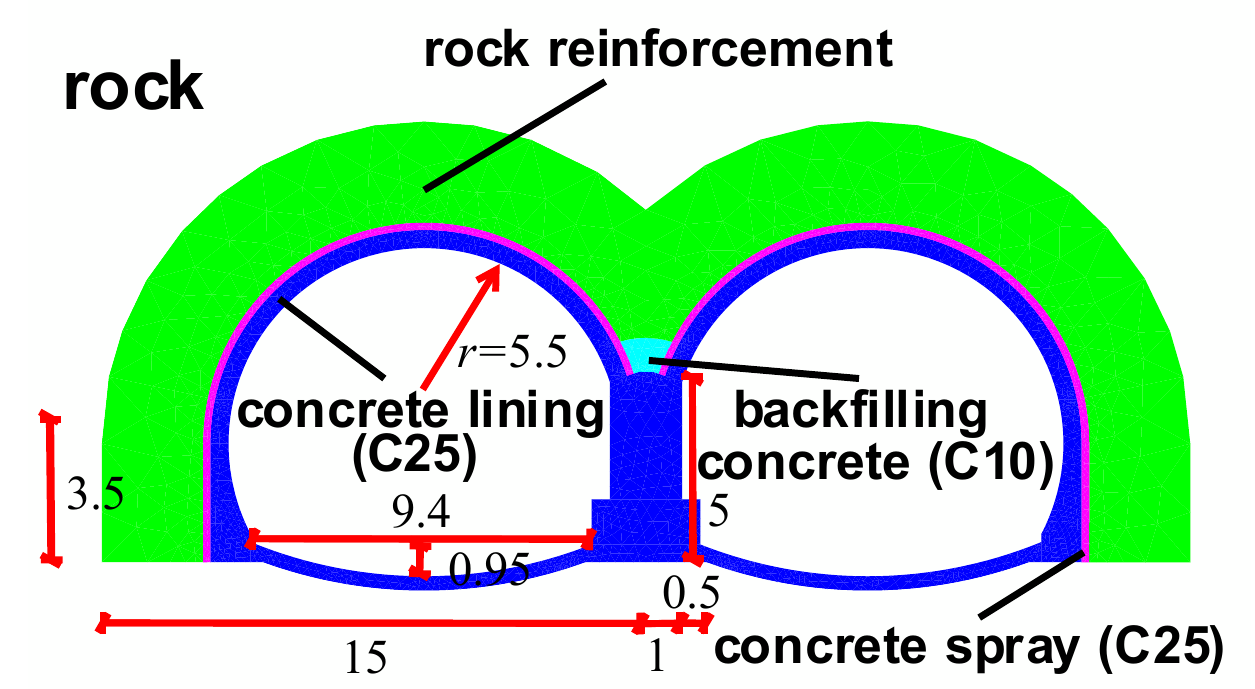}}
	\centerline{\includegraphics[width=0.55\linewidth]{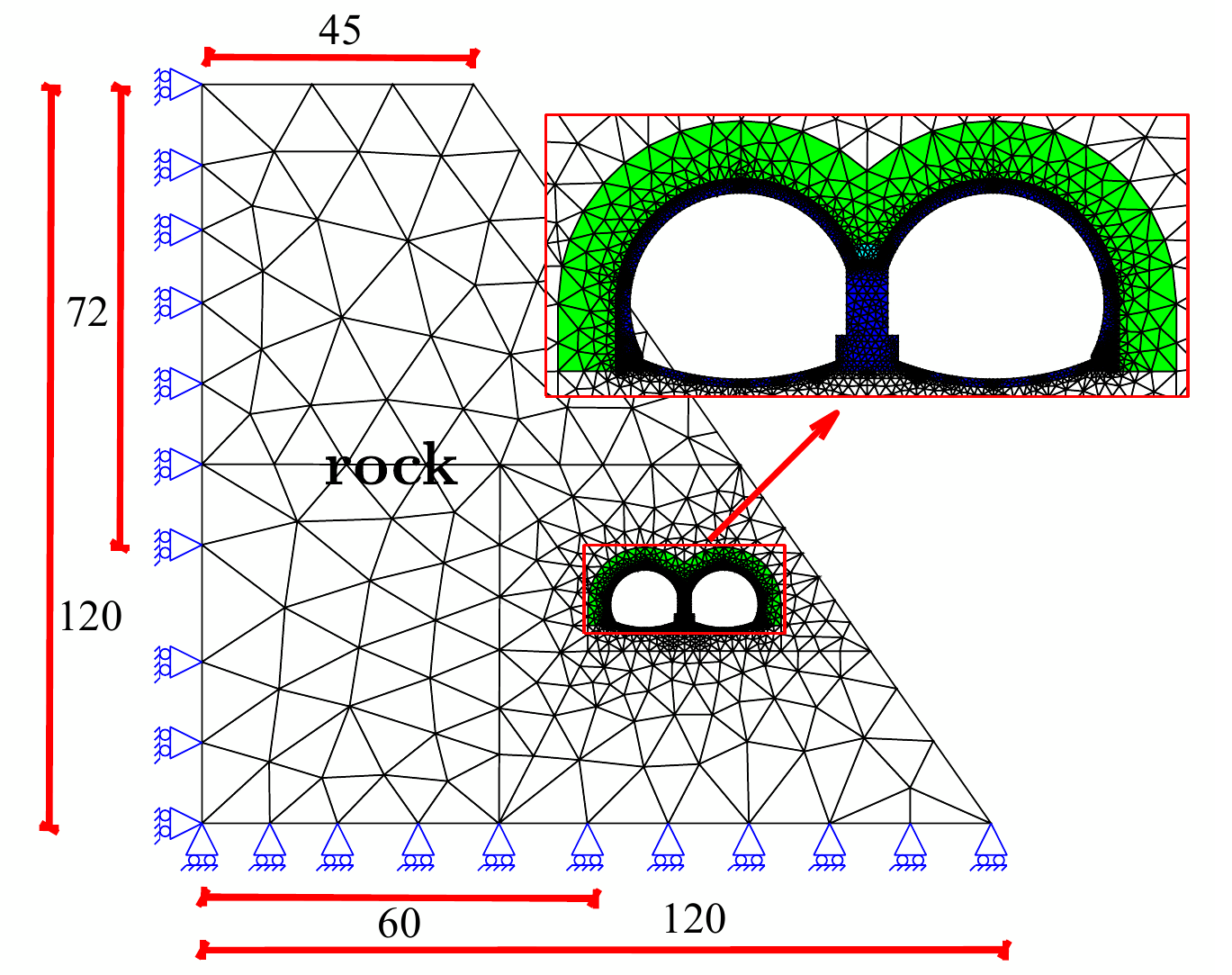}}
	\caption{Model of the tunnel (top) and the finite element mesh (bottom) \label{fig_e3_01}}
\end{figure}

This example is to compute the deformation of a tunnel under the action of self-weight~\citep{XRQ2006}. In order to reduce the size of the stochastic finite element equation while ensuring the accuracy, triangle elements with gradients are used to generate a fine mesh for the tunnel structure and a coarse mesh for the rock, totally including 2729 nodes and 5145 triangle elements, as shown in Figure \ref{fig_e3_01}. Material properties and thicknesses of all components are seen from Table \ref{tab_e3}, here we consider the Young's modulus of components as a random field with the mean value shown in Table \ref{tab_e3} and the covariance function ${C_{EE}}\left( {{x_1},{y_1};{x_2},{y_2}} \right) = \sigma _E^2{e^{ - {{\left| {{x_1} - {x_2}} \right|} \mathord{\left/ {\vphantom {{\left| {{x_1} - {x_2}} \right|} {{l_x}}}} \right. \kern-\nulldelimiterspace} {{l_x}}} - {{\left| {{y_1} - {y_2}} \right|} \mathord{\left/ {\vphantom {{\left| {{y_1} - {y_2}} \right|} {{l_y}}}} \right. \kern-\nulldelimiterspace} {{l_y}}}}}$, where variance function $\sigma _E = 0.1$, correlation lengths ${l_x} = 10$, ${l_y} = 20$. Similar to Example \ref{Example-2}, we model the Young's modulus random field by use of Karhunen-Lo\`{e}ve expansion with 10 terms, and derive a stochastic finite element equation.

\begin{table*}
	\caption{Descriptions of materials properties \label{tab_e3}}
	\centering
	\begin{tabular}{lrrrr}
		\toprule
		& Young's modulus    & Poisson's ratio & mass density  & thickness \\
		&    (GPa)     &       &  ($\rm{{{kg} \mathord{\left/ {\vphantom {{kg} {{m^3}}}} \right. \kern-\nulldelimiterspace} {{m^3}}}}$)  & ($\rm{m}$) \\
		\midrule
		rock &   2.0 & 0.25 & 2200 &  \\
		rock reinforcement &   2.6 & 0.20 & 2300 & 2.80 \\
		concrete lining & 28.5 & 0.20 & 2500 & 0.20 \\
		backfilling concrete & 18.5 & 0.20 & 2300 & 0.50 \\
		concrete spray & 28.5 & 0.20 & 2200 & 0.95 \\
		\bottomrule
	\end{tabular}
\end{table*}

\begin{figure*}
	\begin{minipage}[t]{0.5\linewidth}
		\centerline{\includegraphics[width=1.0\textwidth]{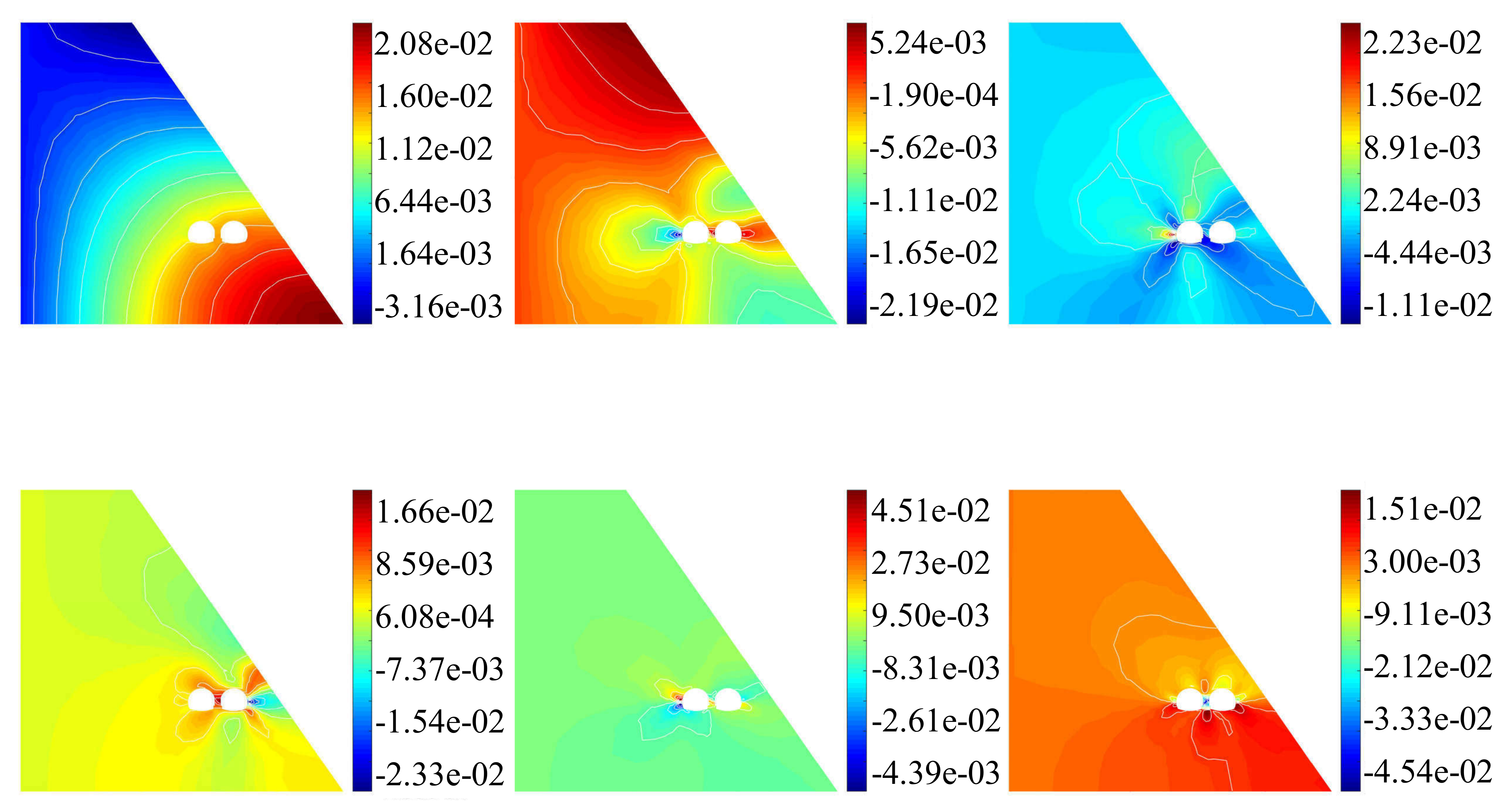}}
		\centerline{$\bf{(a)}$.~Displacement $\left\{ {{d_i}} \right\}_{i = 1}^6$ in $x$ direction \label{fig3_dx}}
	\end{minipage}
	\begin{minipage}[t]{0.5\linewidth}
		\centerline{\includegraphics[width=1.0\textwidth]{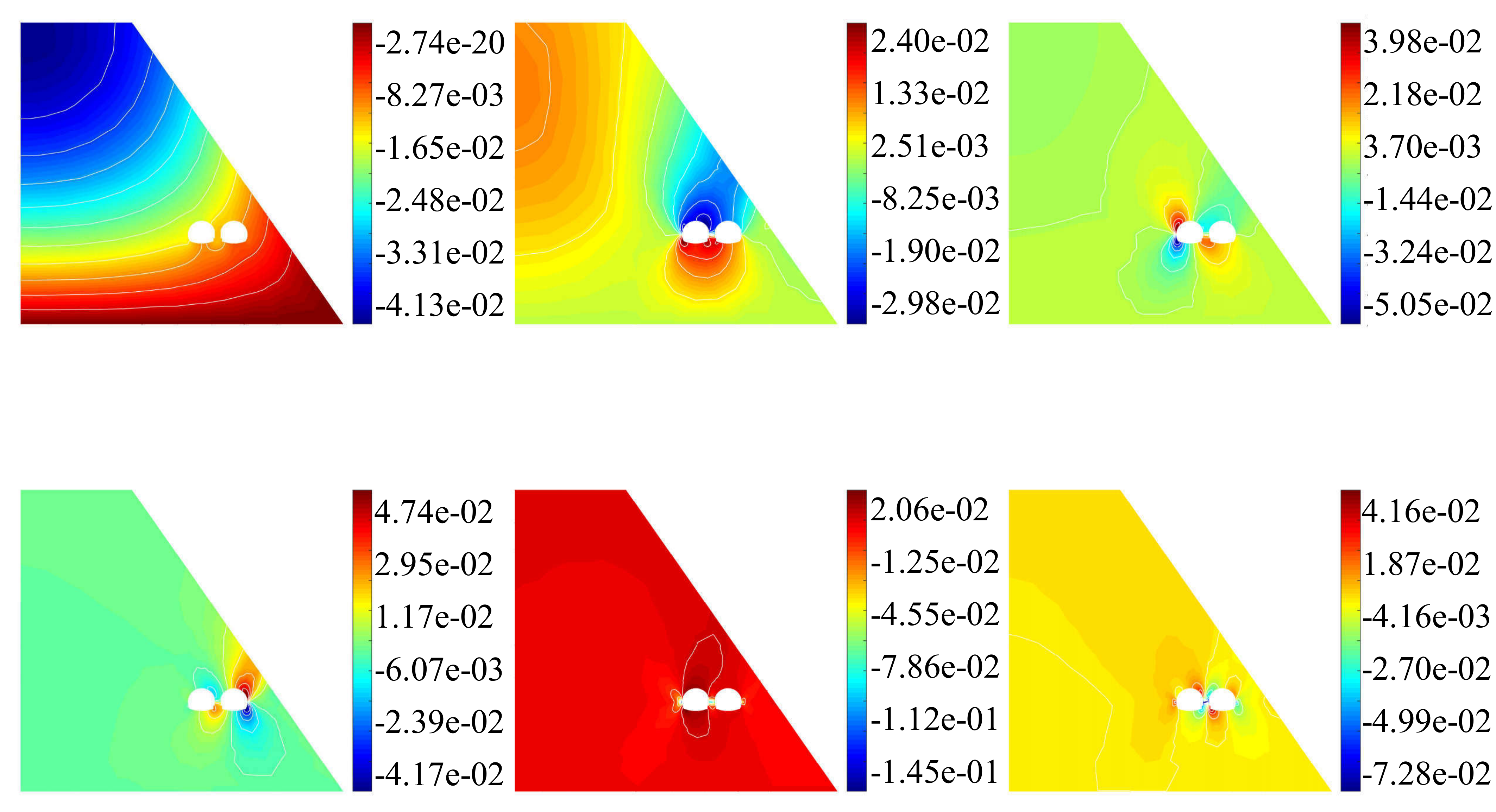}}
		\centerline{~~~~~~$\bf{(b)}$.~Displacement $\left\{ {{d_i}} \right\}_{i = 1}^6$ in $y$ direction \label{fig_3_dy}}
	\end{minipage}
	\begin{minipage}[t]{0.595\linewidth}
		\centerline{\includegraphics[width=1.0\textwidth]{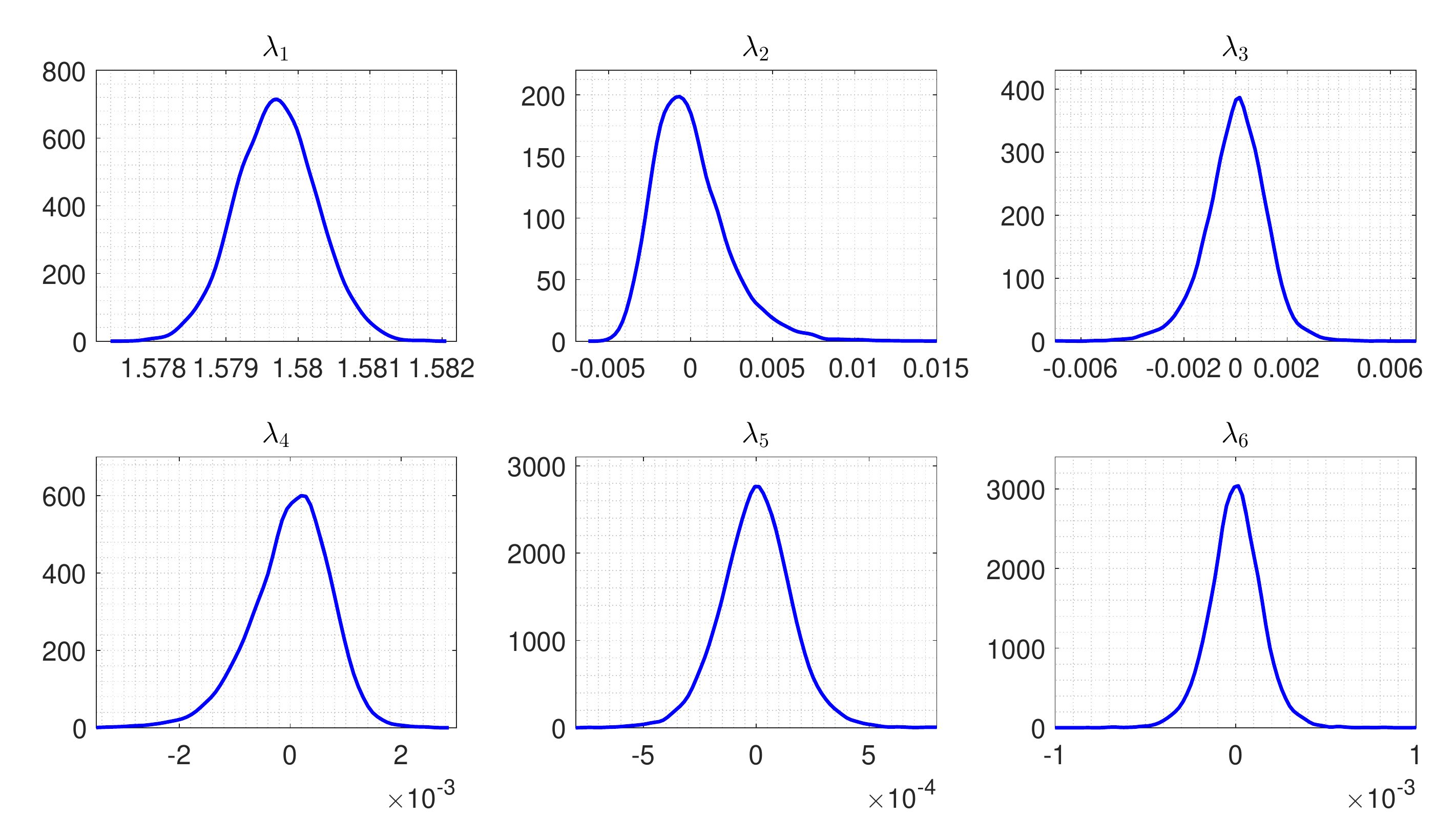}}
		\centerline{$\bf{(c)}$.~PDFs of $\left\{ {{\lambda _i}\left( \theta  \right)} \right\}_{i = 1}^6$ \label{fig3_lam}}
	\end{minipage}
	\begin{minipage}[t]{0.44\linewidth}
		\centerline{\includegraphics[width=1.0\textwidth]{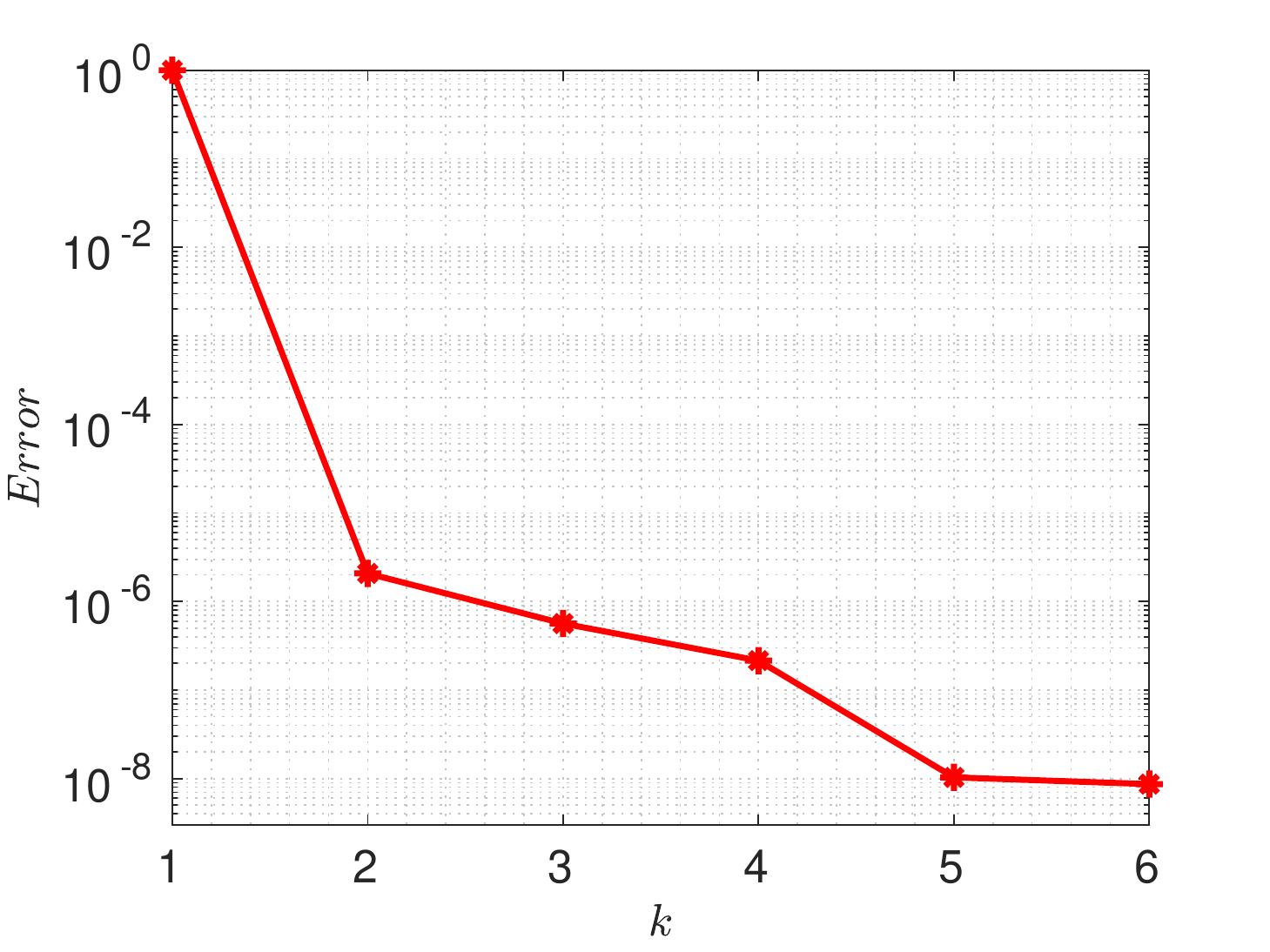}}
		\centerline{$\bf{(d)}$.~Iterative errors of $k$ retained items \label{fig1}}
	\end{minipage}
	\caption{Solutions of the couples $\left\{ {{\lambda _i}\left( \theta  \right),{d_i}} \right\}_{i = 1}^6$ and iterative errors of the solving process \label{fig_e3_02}}
\end{figure*}

\begin{figure*}[!h]
	\begin{minipage}[t]{0.5\linewidth}
		\centerline{\includegraphics[width=1.0\textwidth]{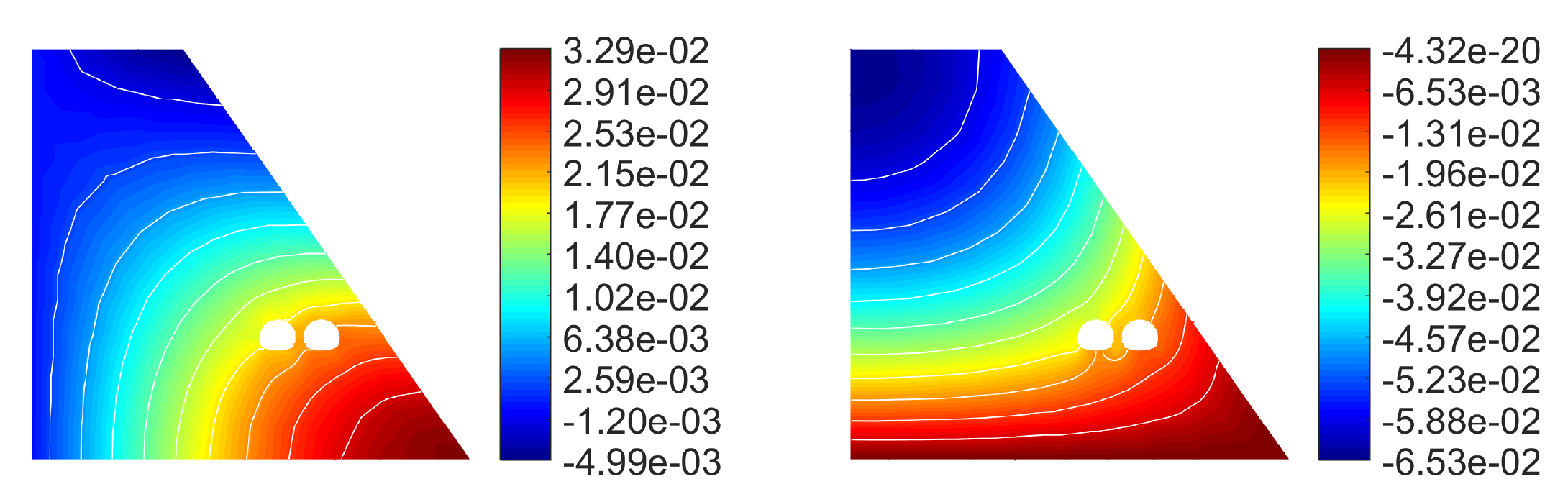}}
		\centerline{$\bf{(a)}$.~Means in the $x$ and $y$ direction \label{fig331}}
	\end{minipage}
	\begin{minipage}[t]{0.5\linewidth}
		\centerline{\includegraphics[width=1.0\textwidth]{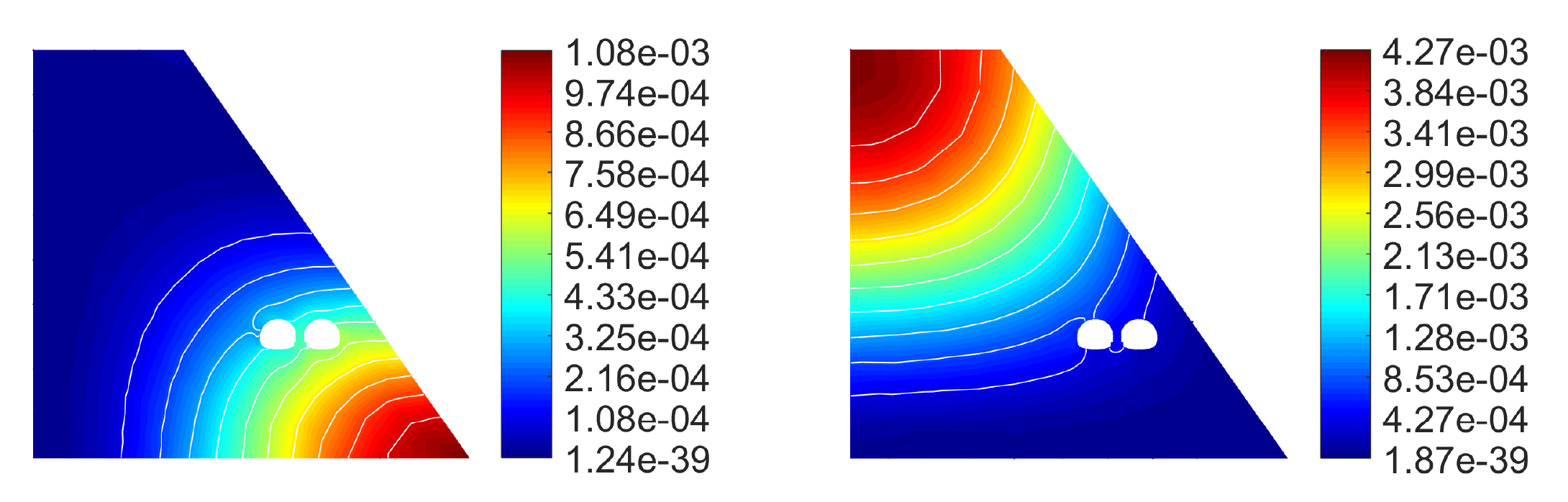}}
		\centerline{$\bf{(b)}$.~Variances in the $x$ and $y$ direction \label{fig332}}
	\end{minipage}
	\begin{minipage}[t]{0.5\linewidth}
		\centerline{\includegraphics[width=1.0\textwidth]{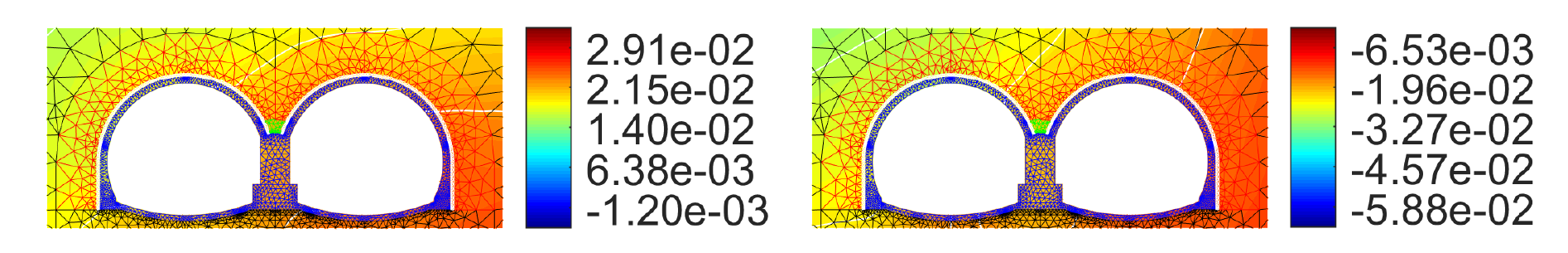}}
		\centerline{$\bf{(a_0)}$.~Means of the tunnel displacement \label{fig333}}
	\end{minipage}
	\begin{minipage}[t]{0.5\linewidth}
		\centerline{\includegraphics[width=1.0\textwidth]{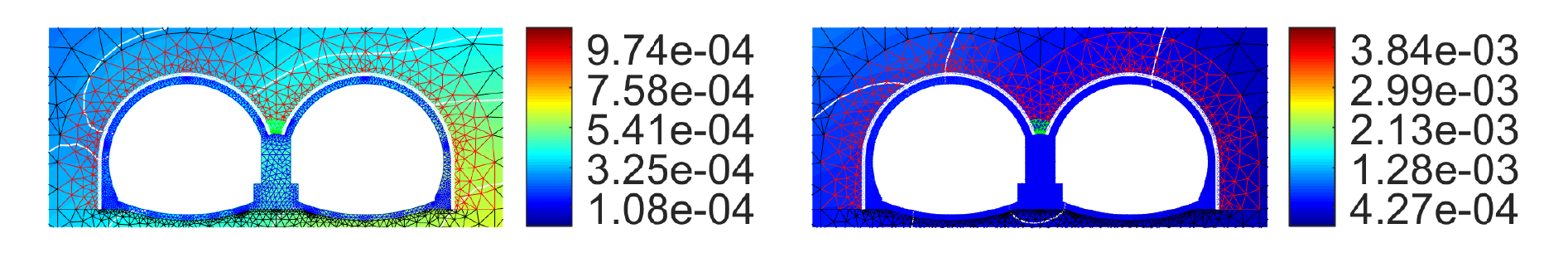}}
		\centerline{~~~~~~~~$\bf{(b_0)}$.~Variances of the tunnel displacement \label{fig334}}
	\end{minipage}
	\caption{Means and variances of the displacement in the $x$ and $y$ direction}
	\label{fig_e3_03}
\end{figure*}

Given the random samples $\left\{ {{\xi _i}\left( {{\theta ^{\left( r \right)}}} \right)} \right\}_{r = 1}^{1 \times {{10}^4}},~i = 1, \cdots ,10$, the initial random variable samples $\{ \lambda _k^{(0)}({\theta ^{(r)}})\} _{r = 1}^{1 \times {{10}^4}}$ and set the convergence criterias as $\varepsilon _1 = 10^{ - 8}$, $\varepsilon _2 = 10^{ - 6}$, displacements $\left\{ {{d_i}} \right\}$ and corresponding random variables $\left\{ {{\lambda _i}\left( \theta  \right)} \right\}$ can be determined after 6 iterations, as shown in Figure \ref{fig_e3_02}{\bf{d}}, which indicates the high efficiency of the proposed method. Figure \ref{fig_e3_02}({\bf{a}}--{\bf{c}}) shows the displacement components $\left\{ {{d_i}} \right\}_{i = 1}^6$ and PDFs of corresponding random variables $\left\{ {{\lambda _i}\left( \theta  \right)} \right\}_{i = 1}^6$, where Figure \ref{fig_e3_02}{\bf{a}} and {\bf{b}} are the displacement components in the $x$ direction (horizontal direction) and $y$ direction (vertical direction), respectively. Mean values and variances of the displacement are shown in Figure \ref{fig_e3_03}$\bf{a}$ and $\bf{b}$, and as a part of the whole displacement (shown in Figure \ref{fig_e3_01} bottom), mean values and variances of the tunnel displacement are seen from Figure \ref{fig_e3_03}$\bf{a_0}$ and $\bf{b_0}$. 
Both tunnel displacements and rock displacements can be captured efficiently, which once again demonstrates the effectiveness of the proposed method. Comparing with Figure \ref{fig_e3_02}($\bf{a}$--$\bf{c}$) and Figure \ref{fig_e3_03}$\bf{a}$, we observe that the first retained item, i.e. $E\left\{ {{\lambda _1}\left( \theta  \right)} \right\}{d_{x1}}$ and $E\left\{ {{\lambda _1}\left( \theta  \right)} \right\}{d_{y1}}$, can roughly approximate the mean displacements in the $x$ direction and $y$ direction. For most cases, the mean displacement considering uncertainties is very close to the displacement obtained from the deterministic case, but considering uncertainties can be better to reflect the variabilities of displacements, which is of great significance for structure design and evaluation, such as reliability analysis and sensitivity analysis~\citep{Dai2014An}. It is seen from Figure \ref{fig_e3_03}$\bf{b}$ that, randomness in this example has more influence on the variance of the tunnel displacement in the $x$ direction and less influence on that in the $y$ direction, which provides a potential way for the design and evaluation of tunnel structures considering uncertainties.

\section{Conclusions}\label{sec6}
In this paper, we develop a method for solving stochastic finite element equations and illustrate its accuracy and efficiency on three practical examples. The proposed method sovles stochastic problems by use of a universal solution construct and a dedicated iterative algorithm. It allows to solve high-dimensional stochastic problems with very low computational costs, which has been illustrated on numerical examples. Thus it appears as a powerful way to avoid the Curse of Dimensionality. In addition, since the stochastic analysis and deterministic analysis in the solving procedure are implemented in their individual spaces, the existing FEM and ODE codes can be readily incorporated into the computatinoal procedure. In these senses, this method is particulary appropriate for large-scale and high-dimensional stochastic problems of practical interests and has great potential in uncertainty quantification of practical problems in science and engineering. In the follow-up research, it hopefully further applies the proposed method to a wider range of uncertainty quantification, such as reliability analysis, sensitivity analysis, etc.

\section*{Acknowledgments}
This research was supported by the Research Foundation of Harbin Institute of Technology and the National Natural Science Foundation of China (Project 11972009). These supports are gratefully acknowledged.

\nocite{*}
\bibliography{References}

\end{document}